\documentclass[12pt]{article}
\usepackage{amsfonts}
\usepackage{amssymb}

\input{tcilatex}
\usepackage{polski}
\begin{document}

\begin{center}
\textbf{Beyond Lebesgue and Baire IV:\\[0pt]
Density topologies and a converse Steinhaus-Weil Theorem\\[0pt]
by \\[0pt]
N. H. Bingham and A. J. Ostaszewski}\\[0pt]
\bigskip

\textit{On the centenary of Hausdorff's Mengenlehre (1914) and Denjoy's
Approximate continuity (1915)}

\bigskip
\end{center}

\noindent \textbf{Abstract.}

The theme here is category-measure duality, in the context of a topological
group. One can often handle the (Baire) category case and the (Lebesgue, or
Haar) measure cases together, by working bi-topologically: switching between
the original topology and a suitable refinement (a density topology). This
prompts a systematic study of such density topologies, and the corresponding
$\sigma $-ideals of negligibles. Such ideas go back to Weil's classic book,
and to Hashimoto's ideal topologies. We make use of \textit{group norms},
which cast light on the interplay between the group and measure structures.\
The Steinhaus-Weil interior-points theorem (`on $AA^{-1}$') plays a crucial
role here; so too does its converse, the Simmons-Mospan theorem.

\bigskip

\noindent \textbf{Key words.} Steinhaus-Weil property, Weil topology,
shift-compact, density topology, Hashimoto ideal topology, group norm

\noindent \textbf{Mathematics Subject Classification (2000): }Primary 26A03;
39B62.

\section{Introduction}

This paper originates from several sources. The first three are our earlier
studies `Beyond Lebesgue and Baire I-III' ([BinO1,2], [Ost3]), the general
theme of which is the similarity (indeed, duality) between measure and
category, and the primacy of category/topology over measure in many areas.
The second three are recent studies by the second author on the Effros Open
Mapping Principle ([Ost4,5,6]; \S 6.4 below).

The \textit{Steinhaus-Weil property} (critical for \textit{regular variation}%
: see [BinGT, Th. 1.1.1]) of a set $S$ in a topological group $G$ is that $%
1_{G}$ is an interior point of $SS^{-1}$ when $S$ is non-negligible (as in
the classic examples in the additive group $\mathbb{R}$: Baire non-meagre,
or measurable non-null). This is implied by the compactness-like property
(called \textit{shift-compactness} in [Ost3], cf. [MilO]), that any \textit{%
null sequence} $z_{n}$ (i.e. $z_{n}\rightarrow 1_{G}$) has a `translator' $%
s\in S$ and subsequence $z_{n(m)}$ with $\{s\cdot z_{n(m)}:m\in N\}\subseteq
S.$ This is not only a stronger property, but also better adapted for use in
many proofs.

Results of Steinhaus-Weil type go back to Steinhaus [Ste] on the line and
Weil [Wei, p. 50] in a locally compact topological group (see e.g.
Grosse-Erdmann [GroE]), and to Kemperman [Kem] (cf [Kuc, Lemma 3.7.2], and
[BinO1, Th. K], where this is `Kemperman's Theorem', [BinO6, Th. 1(iv)]).
For present purposes it is natural to call shift-compactness\ a strong
\textit{Steinhaus-Weil-like} property. A first study of the closeness of the
two Steinhaus-Weil-like properties appears most recently in [Ost6] by way of
the Effros Opening Mapping Principle. Recently, in [BinO8], in work on
subadditivity and mid-point convexity, with `negligibility' interpreted via
the Christensen Haar-null subsets of a Banach space, we replaced
shift-compactness by the already established Steinhaus-Weil property due to
Christensen [Chr1,2], or its extension due to Solecki [Sol2]. (Boundedness
of a subadditive function on $A$ and $B$ yields its boundedness on $AB$ and
hence on an open set, provided $AB$ has the interior-point property -- see
\S 6.9.) This new perspective motivates the present return to the question
of when the Steinhaus-Weil property implies shift-compactness, and hinges on
two themes. The first is the \textit{Lebesgue density theorem}, which
Kemperman [Kem] used to reprove the Steinhaus-Weil theorem, very much a
local property. The second was our reliance [BinO8] on a localized version
of the Steinhaus-Weil property: in $S$ the relative open neighbourhoods of
all points were to have the Steinhaus-Weil property. Taken together, these
suggested the need to \textit{characterize} in a topological group those
analogues of the \textit{Lebesgue density topology} ([HauP], [GofNN],
[GofW]) that imply shift-compactness, a matter we turn to in \S 2. Here we
prove Theorems 1 and 2: we study category analogues of the Lebesgue density
topology. We turn in \S 3 to Hashimoto ideal topologies and in \S 4 to
properties of Steinhaus-Weil type and converses (cf. Prop. 2). We establish
in \S 5 the equivalence of weak and strong Steinhaus-Weil-like properties in
the presence of three topological restrictions, taking what we call the
\textit{Kemperman property} in \S 2 as a weak form of the Steinhaus-Weil
property. This is reminiscent of the characterization of Borel measures on $%
\mathbb{R}$ having the Steinhaus-Weil property, for which see [Mos] ([Sim]
for the Haar case) and the recent [Dan]. We close in \S 6 with some
complements.

We remind the reader of the tension between two of our themes here: density
topologies and topological groups. The real line is not a topological group
under the (Lebesgue) density topology (see e.g. [Sch, Prop. 1.9], [BinO7,
\S\ 4]). Instead, what is relevant here is semi-topological group structure
[ArhT], in which it is the shift (one argument), rather than multiplication
(two arguments) which is continuous.

The Lebesgue Density Theorem, which underlies the density topology (or
topologies) crucial here, is already relevant to (though less well known
than) the Lebesgue Differentiation Theorem, on the relationship between
differentiation and the Lebesgue integral -- see Bruckner's classic survey
[Bru]. The first is usually obtained from the second by specializing to
indicator functions $1_{A}.$ Relevant here are Vitali's covering lemma and
weak \ $L_{1}$-estimates for the Hardy-Littlewood maximal function; for
textbook treatment see e.g. [Rud, Th. 8.8], [SteiS, 3.1.2]. Latent here is
the relation between a general measure $\mu $ and its translation $_{x}\mu $
(defined via $_{x}\mu (B)=\mu (xB)$). That, in turn, is encapsulated in the
Radon-Nikodym derivative $d_{x}\mu /d\mu $ (wherever defined) and is related
to the \textit{relative interior-point property}, which arises when studing
the difference set $S-S$ relative to the \textit{Cameron-Martin subspace} of
a topological vector space; see [Bog, \S 2.2, 2.4]. This is a matter we hope
to return to elsewhere.

In sum: as well as the historical references to Hausdorff in 1914 [Hau] and
Denjoy in 1915 [Den], the paper relates to Lebesgue's approach to the
fundamental theorem of calculus. Its roots may thus be traced back to the
roots of calculus itself.

\section{Density topologies}

Let $(G,\mathcal{T})$ be a separable topological group metrized by a \textit{%
right}-invariant metric $d=d_{R}^{G}$ fixed throughout, allowing the
topology to be denoted either as $\mathcal{T}_{d}$, or $\mathcal{T}$. (This
is possible by the \textit{Birkhoff-Kakutani metrization theorem}, [HewR,
II.8, p. 70 Th., p. 83 Notes], [Ost3].) We make free use of%
\[
||t||:=d_{R}^{G}(1_{G},t),
\]%
referred to as the \textit{group norm} of $G$, for which see the textbook
account in [ArhT, \S 3.3] or [BinO4], and denote by $B_{\delta
}(g):=\{h:||hg^{-1}||<\delta \}=B_{\delta }(1_{G})g$ the open $\delta $-ball
centred at $g$, briefly the open $\delta $-neighbourhood ($\delta $\textit{%
-nhd}) of $g$; we use $B_{\delta }$ to denote $B_{\delta }(1_{G}).$ For $G$
locally compact, we denote (left) Haar measure by $\eta $, or context
permitting by $|.|,$ by analogy with the group norm in view of their close
relationship (cf. \S 6.1). We have in mind, as canonical examples, $\mathbb{R%
}$ or $\mathbb{R}_{+}$ under the usual (Euclidean) topology, denoted $%
\mathcal{E}$. For \textit{any} topology $\mathcal{\tau }$ on $G,$ we write $%
\mathcal{F}(\mathcal{\tau }),\mathcal{F}_{\sigma }(\mathcal{\tau }),\mathcal{%
G}_{\delta }\mathcal{(\tau )}$ for the corresponding closed sets etc., $%
\mathcal{B}(\mathcal{\tau })$ for the \textit{Baire} sets, i.e. the sets
with the Baire Property (BP), $\mathcal{B}_{0}(\mathcal{\tau })$ for the
corresponding meagre sets, and $\mathcal{B}_{+}(\mathcal{\tau })$ for the
non-meagre members of $\mathcal{B}(\mathcal{\tau })$. If $(G,\mathcal{T}$)
is suppressed, $(\mathbb{R}$, $\mathcal{E)}$ or $(\mathbb{R}_{+}$, $\mathcal{%
E)}$ is to be understood. Thus $\mathcal{B}$ denotes the usual Baire sets
and $\mathcal{B}_{0}$ its negligible sets, the $\sigma $-ideal of meagre
sets; analogously, $\mathcal{L}$ denotes the Lebesgue (Haar) measurable sets
and $\mathcal{L}_{0}$ its negligible sets, the $\sigma $-ideal of null
(measure-zero) sets. We denote by $\mathcal{M}(G)$ the Borel regular $\sigma
$-finite measures on $G,$ with $\mathcal{P}(G)$ the subfamily of probability
measures; here regularity is taken to imply both inner and outer regularity
(i.e. compact inner approximation and open outer approximation); these play
a significant role in \S 4. We say that a property holds at \textit{quasi
all points} of a set if it holds except on a negligible set (in the category
or measure sense).

We will refer to the action of $G$ on itself by $t(x)\mapsto tx$ (or $%
t(x)=t+x$ in the case of $\mathbb{R}$ -- we will feel free to move at will
between $(\mathbb{R}$, $\mathbb{+)}$ and $(\mathbb{R}_{+}$, $\cdot )$ via
the exponential isomorphism). Say that a topology $\mathcal{\tau }$ on $G$
is \textit{(left) shift-invariant} if $tV\in \mathcal{\tau }$ for all $t\in
G $ and all $V\in \mathcal{\tau }$; equivalently: each shift $t:\mathcal{%
\tau \rightarrow \tau }$ is continuous.

A \textit{weak} $\mathcal{\tau }$\textit{-base} for $\mathcal{\tau }$ is a
subfamily $\mathcal{W}$ such that for each non-empty $V\in \mathcal{\tau }$
there is $W\in \mathcal{W}$ with $\emptyset \neq W\subseteq V.$ When $%
\mathcal{W}$ above consists of sets \textit{analytic} under $\mathcal{T}_{d}$
(for which see below), the topology is called in [Ost1] a \textit{%
generalized Gandy-Harrington} topology, by analogy with its classical
antecedent (for a textbook treatment of which see [Gao, Ch. 1]); in such a
case the topology $\mathcal{\tau }$ satisfies the Baire Theorem (see [Oxt,
Ch. 9], [Kec, III.26.18,19, p. 203-4], [Ost1, \S 2.2]). Here we consider a
stronger property generalizing the two observations that

\noindent (a): modulo $\mathcal{L}_{0}$ each measurable set is an $\mathcal{F%
}_{\sigma };$

\noindent (b): modulo $\mathcal{B}_{0}$ each Baire set is a $\mathcal{G}_{%
\mathcal{\delta }}$ ([Oxt, Th. 4.4], cf. [Kec, 8.23])

\noindent (these are the forms in which the results are usually stated; it
is the similarities, rather than the distinctions, between $\mathcal{F}%
_{\sigma }$ and $\mathcal{G}_{\delta }$ that are relevant here).

Say that $\mathcal{\tau }$ has the \textit{strong} \textit{Gandy-Harrington}
property if modulo $\mathcal{B}_{0}(\mathcal{\tau })$ each $\mathcal{B}(%
\mathcal{\tau })$ set is analytic under $\mathcal{T}_{d}.$ (Again see the
references above.)

Denote by $\mathcal{D}_{\mathcal{L}}$ the family of all sets $M$ all of
whose points are \textit{density points} (i.e. have Lebesgue (Haar) density
1, in the sense of Martin [Mar1,2] or in the more general context Mueller
[Mue] -- see the more recent development in [Ost3]; cf. [BinO4] for normed
groups, [Oxt2] for $\mathbb{R}$). As noted by Haupt and Pauc [HauP] in $%
\mathbb{R}$, $\mathcal{D}_{\mathcal{L}}$ forms a topology, the (Lebesgue)
density topology. It is related to Denjoy approximate continuity. It can be
generalized to Haar measure. It is a \textit{fine topology} (refining
topology); see [CieLO], [EvaG], [KanK], and [LukMZ], for background on such
fine topologies. (For other topologies derived from notions of `density
point' see [Wil] and e.g. [FilW]; for aspects of translation invariance see
[WilK].)

We list below a number of \textit{qualitative} properties of $\mathcal{D}_{%
\mathcal{L}}$, (i)-(viii), all of them classical. We name property (iv) the
\textit{Kemperman property}: see \S 1 for the motivation, and (vii) the
\textit{Nikodym property }([Ost6]; [Rog, \S 2.9], [Nik]). Property (viii)
suggests a category analogue $\mathcal{D}_{\mathcal{B}}$ of $\mathcal{D}_{%
\mathcal{L}}$; we prove the category analogues of (i)-(vii) in Theorems 1
and 2 below.

\bigskip

\noindent (i) $\mathcal{D}_{\mathcal{L}}$ is a \textit{fine} topology (i.e.
refining $\mathcal{T}_{d}$ ) which is shift-invariant:%
\begin{eqnarray*}
\mathcal{S} &\subseteq &\mathcal{D}_{\mathcal{L}}\Longrightarrow \bigcup
\mathcal{S}\in \mathcal{D}_{\mathcal{L}}, \\
V,V^{\prime } &\in &\mathcal{D}_{\mathcal{L}}\Longrightarrow V\cap V^{\prime
}\in \mathcal{D}_{\mathcal{L}}, \\
\mathcal{T} &\subseteq &\mathcal{D}_{\mathcal{L}}, \\
tV &\in &\mathcal{D}_{\mathcal{L}}\qquad (t\in G,V\in \mathcal{D}_{\mathcal{L%
}});
\end{eqnarray*}%
\noindent (ii) the sets in $\mathcal{D}_{\mathcal{L}}$ are measurable: $%
\mathcal{D}_{\mathcal{L}}\subseteq \mathcal{L}$;

\noindent (iii) the $\mathcal{D}_{\mathcal{L}}$-boundary of a measurable set
is null:%
\[
M\backslash \mathrm{int}_{\mathcal{D}_{\mathcal{L}}}(M)\in \mathcal{L}%
_{0}\qquad (M\in \mathcal{L});
\]%
\noindent (iv) the \textit{Kemperman property,} that any $\mathcal{D}_{%
\mathcal{L}}$-open neighbourhood of the identity meets its own small
displacements non-meagerly: for $1_{G}\in U\in \mathcal{D}_{\mathcal{L}}$
there is $\delta >0$ with
\[
U\cap (tU)\in \mathcal{B}_{+}(\mathcal{D}_{\mathcal{L}})\qquad (||t||<\delta
);
\]

\noindent (v) $\mathcal{D}_{\mathcal{L}}$ is a \textit{strong generalized
Gandy-Harrington topology}: modulo $\mathcal{B}_{0}(\mathcal{D}_{\mathcal{L}%
})$ each $\mathcal{B}(\mathcal{D}_{\mathcal{L}})$ set is analytic under $%
\mathcal{T}_{d}$ and in fact $\mathcal{F}_{\sigma }(\mathcal{T}_{d})\cap
\mathcal{D}_{\mathcal{L}}$ is a weak $\mathcal{D}_{\mathcal{L}}$-base, so
that $\mathcal{D}_{\mathcal{L}}$ is a Baire space;

(Proof: Any set $M\in \mathcal{D}_{\mathcal{L}}$ contains a $\mathcal{T}_{d}$%
-closed set $F$ of positive measure. Let $H$ be a $\mathcal{G(\mathcal{T}}%
_{d})_{\delta }$-null set covering the null set of non-density points$;$
then $F\backslash H\in \mathcal{D}_{\mathcal{L}}\cap \mathcal{F}_{\sigma }(%
\mathcal{T}_{d})$.)

\noindent (vi) the $\mathcal{D}_{\mathcal{L}}$-Baire sets/meagre sets are
identical with respectively, the measurable sets and the null sets:
\[
\mathcal{B}(\mathcal{D}_{\mathcal{L}})=\mathcal{L},\mathcal{\qquad B}_{0}(%
\mathcal{D}_{\mathcal{L}})=\mathcal{L}_{0};
\]

\noindent (vii) the \textit{Nikodym property} of \textit{preservation of
category} under displacements (see [Ost2,6] for background and references):%
\newline
\noindent (a) $tU\in \mathcal{B}(\mathcal{D}_{\mathcal{L}})\qquad (t\in
G,U\in \mathcal{B}(\mathcal{D}_{\mathcal{L}})),$ and \newline
\noindent (b) $tU\in \mathcal{B}_{0}(\mathcal{D}_{\mathcal{L}})$ iff $U\in
\mathcal{B}_{0}(\mathcal{B}_{\mathcal{L}})\qquad (t\in G);$

\noindent (viii) $x$ is a density point of a measurable set $M$ iff $x\in $%
\textrm{int}$_{\mathcal{D}_{\mathcal{L}}}(M).$

So $\mathcal{D}_{\mathcal{L}}$ yields a topological characterization of
local behaviour w.r.t. measure. Property (viii) calls for the Haar
generalization of the Lebesgue Density Theorem [Oxt, Ch. 3]. Below this will
be viewed as a special case of the Banach Localization Principle (or Banach
Category Theorem, [Oxt, Ch. 16]; cf. [Ost1]).

\bigskip

We now define a topology $\mathcal{D}_{\mathcal{B}}$ with properties
analogous to (i)-(viii) in respect of $\mathcal{B}$; here (viii) is a
definition of category-density point.

\bigskip

\noindent \textbf{Definitions. }1. Call $H$ $\mathcal{\tau }$\textit{%
-locally comeagre at }$x\in H$\textit{\ }if there is a $\mathcal{\tau }$%
-open nhd $U$ of $x$ such that $U\backslash H$ is meagre.\textit{\newline
}\noindent 2. Say that $H$ is $\mathcal{\tau }$\textit{-locally comeagre at
all of its points} if for each $x\in H$ there is a $\mathcal{\tau }$-open
nhd $U$ of $x$ such that $U\backslash H$ is meagre.

\bigskip

\noindent \textbf{Remarks. }If, as in (1), $U$ witnesses the property of $H$
at some point $x\in H$, then each point of $H\cap U$ has this property; note
the \textit{monotonicity}: if $H\subseteq H^{\prime }$ and $H$ has the
property at $x,$ then also $H^{\prime }$ has it at $x.$

If (2) holds for $H$, then $H\ $is open under the refinement topology
generated by the family $\{U\backslash L:U\in \mathcal{\tau },L\in \mathcal{B%
}_{0}(\mathcal{\tau })\}.$ We consider `ideal topology' refinements such as
this, generated by a general $\sigma $-ideal in place of $\mathcal{B}_{0}(%
\mathcal{T})$, in the next section.

\bigskip

\noindent \textbf{Theorem 1.} \textit{For }$\mathcal{T}=\mathcal{T}_{d}$%
\textit{\ let }$\mathcal{D}_{\mathcal{B}}(\mathcal{T})$\textit{\ be the
family of sets which are }$\mathcal{T}$\textit{-locally co-meagre at all of
their points. Then:}

\noindent (i) $t\mathcal{D}_{\mathcal{B}}(\mathcal{T})\subseteq \mathcal{D}_{%
\mathcal{B}}(\mathcal{T}),$ \textit{for all} $t\in G$ \textit{and} $\mathcal{%
T}\subseteq \mathcal{D}_{\mathcal{B}}(\mathcal{T});$\textit{\newline
}\noindent (ii)\textit{\ }$\mathcal{D}_{\mathcal{B}}(\mathcal{T})\subseteq
\mathcal{B}(\mathcal{T});$\textit{\newline
}\noindent (iii) $H\backslash \mathrm{int}_{\mathcal{D}_{\mathcal{B}}}(H)\in
\mathcal{B}_{0}(\mathcal{T})$ for $H\in \mathcal{B}(\mathcal{T});$ \textit{%
in fact, if }$S$\textit{\ is Baire, then }$S$\textit{\ has a non-empty }$%
\mathcal{D}_{\mathcal{B}}(\mathcal{T})$\textit{-interior in each nhd on
which it is dense and non-meagre;\newline
}\noindent (iv) \textit{Kemperman property: for} $1_{G}\in H\in \mathcal{D}_{%
\mathcal{B}}(\mathcal{T})$ \textit{there is }$\delta >0$\textit{\ with}
\[
H\cap (tH)\in \mathcal{B}_{+}(\mathcal{D}_{\mathcal{B}}(\mathcal{T}))\qquad
(||t||<\delta );
\]%
\textit{\newline
}\noindent (v) $\mathcal{D}_{\mathcal{B}}$ \textit{is a strong generalized
Gandy-Harrington topology: modulo} $\mathcal{B}_{0}(\mathcal{D}_{\mathcal{L}%
}(\mathcal{T}))$ \textit{each} $\mathcal{B}(\mathcal{D}_{\mathcal{L}}(%
\mathcal{T}))$ \textit{set is analytic under} $\mathcal{T}$\textit{,}
\textit{and in fact} $\mathcal{D}_{\mathcal{B}}(\mathcal{T})\cap \mathcal{G}%
_{\delta }\mathcal{(\mathcal{T})}$\textit{\ is a weak }$\mathcal{D}_{%
\mathcal{B}}(\mathcal{T})$\textit{-base, so }$(G,\mathcal{D}_{\mathcal{B}}(%
\mathcal{T}))$\textit{\ is a Baire space.}

\bigskip

In words: parts (i)-(iii) assert that $\mathcal{D}_{\mathcal{B}}(\mathcal{T}%
) $\textit{\ }is a shift-invariant topology\ refining\textit{\ }$\mathcal{T}$%
\textit{, }the sets in\textit{\ }$\mathcal{D}_{\mathcal{B}}(\mathcal{T})$%
\textit{\ }are\textit{\ }$\mathcal{T}$\textit{-}Baire\textit{, }and\textit{\
}the boundary points of any\textit{\ }$\mathcal{T}$\textit{-}Baire set $H$
form a $\mathcal{T}$\textit{-}meagre set\textit{.}

\bigskip

\noindent \textbf{Proof.} (i) Evidently $G$ and $\emptyset $ are in $%
\mathcal{D}_{\mathcal{B}}(\mathcal{T}).$ For an arbitrary $\mathcal{H}%
\subseteq \mathcal{D}_{\mathcal{B}}(\mathcal{T})$, the union $H:=\bigcup
\mathcal{H}$ is $\mathcal{T}$-locally co-meagre at each element of $H$, by
monotonicity. Next, suppose $x\in H\cap H^{\prime },$ with $H,H^{\prime }\in
\mathcal{D}_{\mathcal{B}}(\mathcal{T})$. Choose $V,V^{\prime }$ open nhds of
$x$ meeting $H,H^{\prime }$ in comeagre sets. Then $x\in W=V\cap V^{\prime }$
is an open nhd of $x$. As $H\cap W$ and $H^{\prime }\cap W$ are co-meagre on
$W,$ by Baire's Theorem, so is their intersection on $W$; so $H\cap
H^{\prime }$ is $\mathcal{T}$-locally co-meagre at $x.$

(ii) Suppose $H$ is $\mathcal{T}$-locally comeagre at all\textit{\ }$x\in H.$%
\textit{\ }Being metrizable and separable, $G$ is hereditarily Lindel\"{o}f
[Dug, VIII \S 6 and Th. 7.3], so $H$ is Baire; indeed, if $\{U_{n}\}$ is a
countable open cover of $H,$ with each $U_{n}\backslash H$ meagre, then $%
H=\bigcup\nolimits_{n}(U_{n}\cap H)=\bigcup\nolimits_{n}U_{n}\backslash
(U_{n}\backslash H),$ which is Baire, since each set $U_{n}\backslash
(U_{n}\backslash H)$ is Baire.

(iii) If $S$ is non-meagre, it is dense in some open nhd $I.$ Then $%
T:=I\backslash S\ $is meagre: for otherwise, $T$ is (Baire and) non-meagre.
Then $T$ is dense on some (non-empty) open $J\subseteq $ $I.$ Being Baire,
both $S$ and $T$ are, modulo meagre sets, $\mathcal{G}_{\delta }(\mathcal{T}%
) $-sets dense on $J.$ So they meet, by Baire's Theorem -- a contradiction.

Let $\mathcal{I}$ be a maximal family of nhds $I$ on which $S$ is co-meagre.
Then $\bigcup \{I\cap S:I\in \mathcal{I\}}$ is $\mathcal{T}$-open. Put $%
S^{^{\prime }}=S\backslash \bigcup \{I\cap S:I\in \mathcal{I\}}$. Then $%
S^{\prime }$ is meagre; otherwise, as before $S^{\prime }$ is dense in some
nhd $I$ and co-meagre on $I,$ contradicting maximality of $\mathcal{I}.$

(iv) If $1_{G}\in H\in \mathcal{D}_{\mathcal{B}}(\mathcal{T}),$ and $%
H=U\backslash N,$ choose $\delta >0$ such that $B_{\delta }(1_{G})\subseteq
U.$ Then $\emptyset \neq \{t\}\subseteq B_{\delta }(1_{G})\cap tB_{\delta
}(1_{G}),$ for $||t||<\delta ,$ so%
\[
H\cap (tH)\supseteq B_{\delta }(1_{G})\cap tB_{\delta }(1_{G})\backslash
(N\cup tN)\in \mathcal{B}_{+}(\mathcal{D}_{\mathcal{B}}(\mathcal{T})).
\]

(v) If $H\in \mathcal{D}_{\mathcal{B}}(\mathcal{T}),$ and $H=U\backslash N$
with $U\in \mathcal{T}$ and $N$ meagre (in the sense of $\mathcal{T}$) we
may chose a larger meagre $\mathcal{F}_{\sigma }(\mathcal{T)}$-set $%
M\supseteq N,$ and then $H^{\prime }=U\backslash M\in \mathcal{G}_{\delta }%
\mathcal{(T)\cap D}_{\mathcal{B}}(\mathcal{T})$. $\square $

\bigskip

Recall the observation of Haupt and Pauc [HauP] that (writing \textit{nwd }%
for nowhere dense)
\begin{equation}
M\text{ is meagre in }\mathcal{D}_{\mathcal{L}}\text{ iff }M\text{ is nwd in
}\mathcal{D}_{\mathcal{L}}\text{ iff }M\text{ is null}  \tag{$H$-$P$}
\end{equation}%
[Kec, 17.47] (the left equivalence as $\mathcal{D}_{\mathcal{L}}$ is a
topology). This shows very clearly how changing from the original topology $%
\mathcal{T}$ ($\mathcal{E}$ in the Euclidean case) to the density topology $%
\mathcal{D}$ turns qualitative measure considerations into Baire (or
topological) considerations. This is the basis for the use of bitopology in
[BinO2].

Theorem 2 below extends the list (i)-(v) in Theorem 1 to (vi), (vii) (as
noted above (viii) becomes a matter of definition), in particular yielding
an abstract form of (H-P).

\bigskip

\noindent \textbf{Theorem 2.} \textit{As in Theorem 1, for }$\mathcal{T}=%
\mathcal{T}_{d}$\textit{:\newline
}\noindent (vi) Haupt-Pauc property:$\ \mathcal{B}_{0}(\mathcal{D}_{\mathcal{%
B}}(\mathcal{T}))=\mathcal{B}_{0}(\mathcal{T})$\textit{\ and }$\mathcal{B}(%
\mathcal{D}_{\mathcal{B}}(\mathcal{T}))=\mathcal{B}(\mathcal{T})$\textit{;%
\newline
}\noindent (vii) Nikodym property:\newline
\noindent (a) $tU\in \mathcal{B}(\mathcal{D}_{\mathcal{B}}(\mathcal{T}))$
\textit{for} $U\in \mathcal{B}(\mathcal{D}_{\mathcal{B}}(\mathcal{T})),$
\textit{and\newline
\noindent }(b) $tU\in \mathcal{B}_{0}(\mathcal{D}_{\mathcal{B}}(\mathcal{T}%
)) $ iff $U\in \mathcal{B}_{0}(\mathcal{B}_{\mathcal{D}}(\mathcal{T}))$.

\bigskip

\noindent \textbf{Proof. }Let $D$ be $\mathcal{D}_{\mathcal{B}}(\mathcal{T})$%
-nwd. Let $\mathcal{W}$ be a maximal family of pairwise disjoint $\mathcal{T}
$-sets $W$ with $D\cap W$ meagre. By Banach's Category Theorem $\bigcup
\{D\cap W:W\in \mathcal{W}\}$ is meagre . Put $D^{\prime }=D\backslash
\bigcup \{D\cap W:W\in \mathcal{W}\}.$ Suppose that $D^{\prime }$ is
non-meagre. Then for some open nhd $U$ the set $D^{\prime }$ is non-meagre
on every open subset of $U;$ but $D^{\prime }$ is $\mathcal{D}_{\mathcal{B}}(%
\mathcal{T})$-nwd, so there is an open nhd $V$ and meagre $M_{V}$ with $%
V\backslash M_{V}\ $disjoint from $D^{\prime }.$ Then $D^{\prime }\cap
V\subseteq M_{V},$ a contradiction. So $D^{\prime }$ is indeed meagre, and
so is $D.$

If $D$ is a countable union of $\mathcal{D}_{\mathcal{B}}(\mathcal{T})$-nwd
sets, then it is meagre in $\mathcal{T}$.

The remaining assertions are now clear. $\square $

\section{Hashimoto ideal topologies}

As in ($H$-$P$) above, a key role is played by the family of meagre sets (in
$\mathcal{E},$ or $\mathcal{B}_{\mathcal{L}}$) and (Lebesgue-) null sets.
Each forms a $\sigma $-ideal of \textit{small sets} and gives rise to a
\textit{Hashimoto topology} (= `ideal topology' according to [LukMZ, 1.C]),
to which we turn in this section. For an illustration of its use in a Banach
space exploiting the $\sigma $-ideal $\mathcal{HN}$ of Haar-null subsets
(defined in \S 5 below, cf. \S 6.6) -- see [BinO8]. Also relevant here is
the study [CieJ] of topologies $\tau $ for which a given $\sigma $-ideal is
identical with the $\sigma $-ideal $B_{0}(\tau )$ of $\tau $-meagre sets.

\bigskip

\noindent \textbf{Definition. }For $\mathcal{I}$ a $\sigma $-ideal of
subsets of $X$ and $\mathcal{\tau }$ a topology on $X,$ say that a set $S$
is $\mathcal{I}$\textit{-nearly }$\mathcal{\tau }$\textit{-open} if for some
$\mathcal{\tau }$-open set $U$ and elements $M,N$ in $\mathcal{I}$%
\[
S=(U\backslash N)\cup M.
\]%
We make the blanket assumption that $\{x\}\in \mathcal{I}$ for all $x\in X$
-- see [Hay] , and also [Sam, especially Ex. 2]).

When relevant, we will include the axioms of set theory we use in
parentheses after the theorem number. As usual, we write ZF for
Zermelo-Fraenkel and DC for Dependent Choice.

The following result is an embellishment of Hashimoto's early insight [Has];
there is a readable account by Jankovi\'{c} and Hamlett in [JanH]
introducing the topology via a Kuratowski closure operator, rather than via
an ideal, as was first done in [Sam], albeit anticipated by [Sch] using $%
\mathcal{L}$ in $\mathbb{R}$, as noted in the Introduction.

\bigskip

\noindent \textbf{Theorem 3} (ZF+DC; cf. [Has], [JanH], [Sam, Cor 2]).
\textit{For a topology }$\mathcal{\tau }$ \textit{with countable basis} $%
\beta $ \textit{and a }$\sigma $\textit{-ideal} $\mathcal{I}$\textit{, the
(Hashimoto) topology generated by the sets }%
\[
\mathcal{H}:=\{V\backslash M:V\in \tau ,M\in \mathcal{I}\}
\]%
\textit{\ is }$\mathcal{H}$\textit{\ itself, so that, in particular, each} $%
W\in \mathcal{H}$ \textit{is }$\mathcal{I}$\textit{-nearly }$\mathcal{\tau }$%
\textit{-open and has the representation}%
\[
W=\left( \bigcup\nolimits_{B\in \mathcal{\beta }_{W}}B\right) \backslash M,%
\text{ for some }M\in \mathcal{I},
\]%
\textit{where}%
\[
\mathcal{\beta }_{W}:=\{B\in \beta :B\backslash L\subseteq W\text{ for some }%
L\in \mathcal{I}\}.
\]%
\textit{Furthermore, if no non-empty }$\mathcal{\tau }$\textit{-open set is
in }$\mathcal{I}$\textit{, then:}

(i)\textit{\ a set is }$\mathcal{H}$\textit{-nowhere dense iff it is the
union }$N\cup M$\textit{\ of a }$\mathcal{\tau }$\textit{-nowhere dense set
and an }$\mathcal{I}$ \textit{set;}

(ii) \textit{the }$\mathcal{H}$\textit{-Baire sets are the }$\mathcal{I}$%
\textit{-nearly }$\mathcal{\tau }$\textit{-Baire sets;}

(iii) \textit{for }$\mathcal{I}$ \textit{the meagre sets of }$\mathcal{\tau }
$\textit{, if }$\mathcal{\tau }$\textit{\ is a Baire topology, then so is }$%
\mathcal{H}$\textit{.}

\bigskip

\noindent \textbf{Proof.} For a fixed non-empty set $W$ that is a union of
sets in $\mathcal{H}$, put%
\[
\mathcal{\beta }_{W}:=\{B\in \mathcal{\beta }:B\backslash L\subseteq W\text{
for some }L\in \mathcal{I}\}.
\]%
As this is countable, by DC\ (see [TomW, Ch. 15]) there is a selector $L(.)$
defined on $\mathcal{\beta }_{W}$ such that $L(B)\in \mathcal{I}$ for each $%
B\in \mathcal{\beta }_{W}$; then $B\backslash L(B)\subseteq W$. Set $%
M(B):=B\cap L(B)\backslash W,$ so that $B\backslash M(B)=B\backslash
L(B)\cup (B\cap L(B)\cap W)\subseteq W$. Then $M(B)\cap W=\emptyset .$ So
putting%
\[
M:=\bigcup\nolimits_{B\in \mathcal{\beta }_{W}}M(B),
\]%
$M\cap W=\emptyset $. Also $M\in \mathit{\ }\mathcal{I}$ as $M(B)\subseteq
L(B)\in \mathit{\ }\mathcal{I}$. Then, as $B\backslash M(B)\subseteq W$ for
each $B\in \mathcal{\beta }_{W},$%
\[
\left( \bigcup\nolimits_{B\in \mathcal{\beta }_{W}}B\right) \backslash
M\subseteq W.
\]%
In fact, equality holds: indeed, if $x\in W,$ then $x\in B\backslash L$ for
some $B\in \mathcal{\beta }_{W}$ and $L\in \mathcal{I}$ (and $B\backslash
L\subseteq W).$ Then $x\notin M$ as $x\in W,$ and so $x\in B\backslash M.$

So in particular $W$ is $\mathcal{I}$-nearly $\mathcal{\tau }$-open.

(i) This follows since the $\mathcal{I}$-nearly $\mathcal{\tau }$-open sets
are a $\sigma $-ideal containing $\mathcal{H}$.

(ii) If $N$ is $\mathcal{H}$-nowhere dense, then so is its $\mathcal{H}$%
-closure $\bar{N}$. So $W:=X\backslash \bar{N}$ is $\mathcal{H}$-open and $%
\mathcal{H}$-everywhere dense. Write $W=(V\backslash M)$ with $V\in \mathcal{%
\tau }$ and $M\in \mathcal{I}$. It is enough to show now that $V$ is $%
\mathcal{\tau }$-everywhere dense. Indeed, for $\mathcal{\tau }$-open $U,$
as $U$ is also $\mathcal{H}$-open, $U\cap (V\backslash M)=(U\cap
V)\backslash M$ is non-empty, as otherwise $U\cap V\subseteq M,$ a
contradiction.

Conversely, if $V$ is $\mathcal{\tau }$-everywhere dense and $M\in \mathcal{I%
}$, then we are to show that $(X\backslash V)\cup M$ is $\mathcal{H}$%
-nowhere dense. For non-empty $\mathcal{\tau }$-open $U$ and $M^{\prime }\in
\mathcal{I}$, the set $(U\backslash M^{\prime })\cap (V\backslash M)=(U\cap
V)\backslash (M\cup M^{\prime })$ is non-empty (otherwise $M\cup M^{\prime }$
covers the non-empty set $U\cap V)$. So $V\backslash M$ is $\mathcal{H}$%
-everywhere dense, and so $X\backslash (V\backslash M)$ is $\mathcal{H}$%
-nowhere dense and $X\backslash (V\backslash M)=(X\backslash V)\cup M.$

(iii) Straightforward, since if $W_{n}$ for $n\in \mathbb{N}$ is $\mathcal{H}
$-everywhere dense, then $X\backslash W_{n}$ is $\mathcal{H}$-nowhere dense,
and so%
\[
X\backslash W_{n}=X\backslash G_{n}\cup M_{n}
\]%
for some $G_{n}$ $\mathcal{\tau }$-open $\mathcal{\tau }$-everywhere dense,
and $M_{n}\in \mathcal{I}.$ That is: the topology is Baire. $\square $

\bigskip

\noindent \textbf{Cautionary example.} For $\mathcal{I=L}_{0},$ the Lebesgue
null sets, decompose $\mathbb{R}$ into a meagre and a null set [Oxt, Th.
1.6] to see that starting with $\tau $ the usual topology of $\mathbb{R}$
yields a Hashimoto topology that is not Baire, by (i) above.

\section{Steinhaus-Weil-like properties}

The Kemperman property (iv) of Theorem 1 above, for the $\mathcal{D}$-open
sets $U$:%
\[
U\cap (tU)\in \mathcal{B}_{+}(\mathcal{D})\qquad (||t||<\delta ),
\]%
yields the Steinhaus-Weil property $B_{\delta }\subseteq A^{-1}A$ for a set $%
A$ (see e.g. [BinO6], [BinGT, 1.1.2], [BinO9]), provided the Kemperman
property extends to $A,$ since%
\[
ta=a^{\prime }\qquad \Leftrightarrow \qquad t=a^{\prime }a^{-1}\in AA^{-1}.
\]%
So we regard the Kemperman property as a weak \textit{Steinhaus-Weil-like}
property. The extension beyond $\mathcal{D}$-sets to $\mathcal{B}_{+}(%
\mathcal{D})$ will indeed hold in the presence of the \textit{Nikodym
property} (\S 2 above), when $U\backslash A\in \mathcal{B}_{0}(\mathcal{D})$
for some $U\in \mathcal{D}$, as then $t(U\backslash A)\in \mathcal{B}_{0}(%
\mathcal{D}).$ The fact that this is so for $A$ a Haar-measurable set, or a
Baire set, is the source of the most direct proofs, recalled below for
completeness, of the Steinhaus-Pettis Theorem (see e.g. [BinO4, Th. 6.5]).
These generalize (from the Euclidean case) a Lemma first observed by
Kemperman ([Kem]; cf. [Kuc, Lemma 3.7.2]).

In Lemma 2 below we take a more direct approach -- a streamlined version of
the `uniformity'\ approach in [Hal, Th. 61.A] -- inspired by Stromberg
[Str], but with \textit{translational} \textit{subcontinuity }of measure $%
\mu \in \mathcal{M}(G)$ (its definition below motivated by the upper
semicontinuity of the map\textit{\ }$x\mapsto \mu (xK)$), in place of
\textit{translation-invariance} of measure, and referring to the group norm,
more thematic here. This is followed by its Baire-category analogue. Below,
for $\eta $ a left Haar measure of a locally compact group $G,$ $\mathcal{M}%
_{+}(\eta )$ denotes the left Haar measurable sets of positive (finite)
measure, by analogy with the notation $\mathcal{B}_{+}=\mathcal{B}_{+}(%
\mathcal{T})$ for the non-meagre Baire sets.

A key tool is provided by a form of the `telescope' or `tube' lemma (cf.
[Mun, Lemma 5.8]). Our usage of upper semicontinuity in relation to
set-valued maps follows [Rog], cf. [Bor].

\bigskip

\noindent \textbf{Lemma 1}. (cf. [BeeV]). \textit{For compact }$K\subseteq
G, $ \textit{the map }$t\mapsto tK$ \textit{is upper semicontinuous; in
particular, for }$\mu \in \mathcal{M}(G)$\textit{,}%
\[
\mathit{\ }m_{K}:t\mapsto \mu (tK)
\]%
\textit{is upper semicontinuous, hence }$\mu $\textit{-measurable. In
particular, if }$m_{K}(t)=0$\textit{, then }$m$ \textit{is continuous at }$t$%
\textit{.}

\bigskip

\noindent \textbf{Proof. }For $K$ compact and $V\supseteq K$ open, pick for
each $k\in K,$ an $r(k)>0$ with $kB_{2r(k)}$ $\subseteq V$ (possible, as $G$
is a topological group, by continuity at $x=1_{G}$ of $x\mapsto kx$). By
compactness, there are $k_{1},...,k_{n}$ with $K\subseteq
\dbigcup\nolimits_{j}B_{r(k_{k})}k_{j}\subseteq V;$ then for $\delta
:=\min_{j}r(k_{j})>0$%
\[
tK\subseteq \dbigcup\nolimits_{j}tB_{r(k_{k})}k_{j}\subseteq
\dbigcup\nolimits_{j}B_{2\cdot r(k_{k})}k_{j}\subseteq V\qquad (||t||<\delta
).
\]

To prove upper semicontinuity of $m_{K},$ fix $t\in G.$ For $\varepsilon >0$%
, as $tK$ is compact, choose by outer regularity an open $U\supseteq tK$
with $\mu (U)<\mu (tK)+\varepsilon ;$ by the first assertion, there is an
open ball $B_{\delta }$ at $1_{G}$ with $B_{\delta }tK\subseteq U$, and then
$\mu (B_{\delta }(t)K)\leq \mu (U)<\mu (K)+\varepsilon .$ $\square $

\bigskip

Of course, in a locally compact group $G$ with Haar measure $\mu $, $m_{K}$
is \textit{constant}. By Luzin's theorem, the restriction of $m_{K}$ (for $%
K\ $compact) to appropriate non-null subsets of $G$ is (relatively)
continuous; but of greater significance, as emerges in [BinO9], is a form of
subcontinuity \textit{relativized to a fixed sequence} $t_{n}\rightarrow
1_{G}$, linking the concept to Solecki's amenability at 1 [Sol2] (see \S %
8.5), and the latter, like \textit{outright continuity at} $1_{G}$, yields
the Steinhaus-Weil property (of non left-Haar null universally measurable
sets). Regularity of measure also plays a part; it is likewise key in
establishing in [Kom] the connection between the (wider) Steinhaus-Weil
property concerning $AB^{-1}$ (for which see \S 6.9) and certain forms of
\textit{metric transitivity }of measure (specifically, the co-negligibility/
`residuality' of $AD$ for any countable dense $D$ cf. [Kuc, Th. 3.6.1],
[CichKW, Ch. 7]), known as the \textit{Sm\'{\i}tal property }([KucS], cf.
[BarFN]). We begin with a

\bigskip

\noindent \textbf{Cautionary example.} In any separable group $G,$ for $%
D:=\{d(n):n\in \mathbb{N}\}$ a dense subset, define a regular probability
measure by $\mu _{D}:=\sum\nolimits_{n\in \mathbb{N}}2^{-n}\delta _{d(n)},$
with $\delta _{g}$ the Dirac measure at $g$ (unit point-mass at $g);$ then $%
\mu _{D}(U)>0$ for all non-empty $U$, as $\mu _{D}(d(n))=2^{-n}.$ However,
there exist arbitrarily small translations $t$ with $\mu _{D}(tD)=0$ (since $%
D$ is meagre -- cf. [MilO]). This is particularly obvious in the case $G=%
\mathbb{R}$ with $D=\mathbb{Q}$, on taking small irrational translations;
while the situation here is attributable to $\mu _{D}$ having atoms (not
just being singular w.r.t. Lebesgue measure -- see the Mospan property
below), its obverse occurs when $\mu (U)=0,$ for some non-atomic $\mu \in
\mathcal{P}(G)$ and $U$ non-empty open, as there is a translation $t=d(n)$
such that $\mu (tU)>0$ (as $\{d(n)U:n\in \mathbb{N}\}$ covers $G).$ Then $%
_{t}\mu (U)>0$ for $_{t}\mu (\cdot )=\mu (t\cdot ),$ which is non-atomic.

Motivated by this last comment, and with Proposition 3 below in mind, let us
note a non-atomic modification $\mu $ of $\mu _{D}$ making all open sets $%
\mu $-non-null. In $\mathcal{P}(G),$ under its (separable, metrizable) weak
topology, the non-atomic measures form a dense $\mathcal{G}_{\delta }$ (see
e.g. [Par, II.8]). So for a non-atomic $\tilde{\mu}$ take $\mu
:=\sum\nolimits_{n\in \mathbb{N}}2^{-n}$ $_{d(n)}\tilde{\mu},$ which is
non-atomic with $\mu (U)>0$ for non-empty open $U.$

It is inevitable that the significance of small changes in measure relates
to amenability (see the Reiter condition in [Pat, Prop. 0.4], cf. \S 6.5).

\bigskip

\noindent \textbf{Definition.} For $\mu \in \mathcal{P}(G)$ and compact $%
K\subseteq G,$ noting that $\mu _{\delta }(K):=$ $\inf \{\mu (tK):t\in
B_{\delta }\}$ is weakly decreasing in $\delta $, put%
\[
\mu _{-}(K):=\sup_{\delta >0}\inf \{\mu (tK):t\in B_{\delta }\}.
\]%
Then, as $1_{G}\in B_{\delta }$,%
\[
0\leq \mu _{-}(K)\leq \mu (K).
\]%
We will say that the measure $\mu $ is \textit{translation-continuous,} or
just \textit{continuous, }if $\mu (K)=\mu _{-}(K)$ for all compact $K;$
evidently $m_{K}(.)$ is continuous if $\mu $ is translation continuous,
since $m_{K}(st)=m_{tK}(s)$ and $tK$ is compact whenever $K$ is compact. For
$G$ locally compact this occurs for $\mu =\eta ,$ the left-Haar measure, and
also for $\mu $ absolutely continuous w.r.t. to $\eta $ (see below). We call
$\mu $ \textit{maximally discontinuous} at $K$ if $0=\mu _{-}(K)<\mu (K).$
That a measure $\mu $ singular w.r.t. Haar measure is just such an example
was first discovered by Simmons [Sim] (and independently, much later, by
Mospan [Mos]). The intermediate situation when $\mu (K)\geq \mu _{-}(K)>0$
for all $\mu $-non-null compact $K$ is of significance; then we call $\mu $%
\textit{\ subcontinuous.} The notion of subcontinuity for functions goes
back to [Ful] (cf. [Bou]): as applied to the function $m_{K}(t),$ regarded
as a map into the positive reals $(0,\infty ),$ subcontinuity at $t=1_{G}$
requires that for \textit{every sequence} $t_{n}\rightarrow 1_{G}$ there is
a subsequence $t_{m(n)}$ with $m_{K}(t_{m(n)})$ convergent (to a positive
value). Thus our usage, applied to a measure $\mu ,$ is equivalent to
demanding, for each $\mu $-non-null compact $K$ and \textit{any} null
sequence $\{t_{n}\}$ (i.e. with $t_{n}\rightarrow 1_{G}),$ that there be a
subsequence $\mu (t_{m(n)}K)$ bounded away from $0.$ This is already
reminiscent of \textit{amenability at} $1_{G}$ -- again see \S 6.5. In Lemma
2 below (H-K) is for `Haar-Kemperman', as in Proposition 1 below.

\bigskip

\noindent \textbf{Lemma 2.} \textit{Let }$\mu \in \mathcal{P}(G)$\textit{.
For }$\mu $\textit{-non-null compact }$K\subseteq G,$ \textit{if }$\mu
_{-}(K)>0$ \textit{(i.e. }$\mu $ \textit{is subcontinuous at }$K$),\textit{\
then there is }$\delta >0$\textit{\ with}%
\begin{equation}
tK\cap K\in \mathcal{M}_{+}(\mu )\qquad (||t||<\delta ),  \tag{$H$-$K$}
\end{equation}%
\textit{so in particular,}%
\[
B_{\delta }\subseteq KK^{-1}
\]%
\textit{(so that }$B$\textit{\ has compact closure)}, \textit{or,
equivalently,\ }%
\[
tK\cap K\neq \emptyset \qquad (||t||<\delta ).
\]%
\noindent \textbf{Proof.} Choose $\Delta >0$ with $\mu _{\Delta }(K)>\mu
_{-}(K)/2.$ By outer regularity of $\mu $, choose $U$ open with $K\subseteq
U $ and $\mu (U)<\mu (K)+\mu _{-}(K)/2.$ By upper semicontinuity of $%
t\mapsto tK$, w.l.o.g. $B_{\delta }K\subseteq U$ for some $\delta <\Delta .$
Then $(H$-$K)$ holds for this $\delta $: otherwise, for some $t\in B_{\delta
},$ as $\mu (tK\cap K)=0,$
\[
\mu _{\delta }(K)+\mu (K)\leq \mu (tK)+\mu (K)=\mu (tK\cup K)\leq \mu
(U)<\mu (K)+\mu _{-}(K)/2,
\]%
so $\mu _{-}(K)/2<\mu _{\Delta }(K)<\mu _{\delta }(K)<\mu _{-}(K)/2,$ a
contradiction. Given $||t||<\delta $ and $tK\cap K\in \mathcal{M}_{+}$, take
$s\in tK\cap K\neq \emptyset ;$ then $s=ta$ for some $a\in K,$ so $%
t=sa^{-1}\in KK^{-1}.$ Conversely, $t\in B_{\delta }\subseteq KK^{-1}$
yields $t=a^{\prime }a^{-1}$for some $a,a^{\prime }\in K;$ then $a^{\prime
}=ta\in K\cap tK$. $\square $

\bigskip

Proposition 1M/1B below, which occurs in two parts (measure and Baire
category cases), unifies and extends various previous results due to among
others Steinhaus, Weil, Kemperman, Kuczma, Stromberg, Weil, Wilczy\'{n}ski,
Simmons [Sim] and Mospan [Mos] -- see [BinO6,9] for references.

\bigskip

\noindent \textbf{Proposition 1M }(\textbf{Haar-Kemperman property).}
\textit{Let }$\mu \in \mathcal{P}(G)$\textit{\ with each map }$%
m_{K}:t\mapsto \mu (tK),$\textit{\ for non-null compact }$K\subseteq G,$
\textit{continuous at }$1_{G}$.

\noindent \textit{Then for\ }$\mu $\textit{-measurable }$A$ \textit{with }$%
0<\mu (A)<\infty $\textit{\ there is }$\delta >0$\textit{\ with}%
\begin{equation}
tA\cap A\in \mathcal{M}_{+}(\mu )\qquad (||t||<\delta ),  \tag{$H$-$K$}
\end{equation}%
\textit{so in particular,}%
\[
B_{\delta }\subseteq AA^{-1},
\]%
\textit{or, equivalently,\ }%
\[
tA\cap A\neq \emptyset \qquad (||t||<\delta ).
\]%
\noindent \textbf{Proof.} By inner regularity of $\mu $, there exists a
compact $K\subseteq A$ with $0<\mu (K)\leq \mu (A)<\infty .$ Now apply the
previous lemma. $\square $

\bigskip

We can now give a general form of a result of Mospan (sharpened by provision
of a converse).

\bigskip

\textbf{\noindent Corollary 1 (Mospan property, }[Mos, Th. 2]\textbf{). }%
\textit{For }$\mu $\textit{-non-null compact }$K,$ \textit{if }$1_{G}\notin
\mathrm{int}(KK^{-1}),$\textit{\ then }$\mu _{-}(K)=0,$\textit{\ i.e. }$\mu $
\textit{is maximally discontinuous at }$K$\textit{; equivalently, there is a
`null sequence' }$t_{n}\rightarrow 1_{G}$\textit{\ with }$\lim_{n}\mu
(t_{n}K)=0.$\textit{\ }

\textit{Conversely, if }$\mu (K)>\mu _{-}(K)=0,$ \textit{then there is a
null sequence }$t_{n}\rightarrow 1_{G}$\textit{\ with }$\lim_{n}\mu
(t_{n}K)=0,$\textit{\ and there is }$C\subseteq K$\textit{\ with }$\mu
(K\backslash C)=0$\textit{\ with }$1_{G}\notin \mathrm{int}(CC^{-1}).$

\bigskip

\noindent \textbf{Proof.} The first assertion follows from Lemma 2. For the
converse, as in [Mos]: suppose that $\mu (t_{n}K)=0,$ for some sequence $%
t_{n}\rightarrow 1_{G}.$ By passing to a subsequence, we may assume that $%
\mu (t_{n}K)<2^{-n-1}.$ Put $D_{m}:=K\backslash \bigcap\nolimits_{n\geq
m}t_{n}K\subseteq K;$ then $\mu (K\backslash D_{m})\leq \sum_{n\geq m}\mu
(t_{n}K)<2^{-m},$ so $\mu (D_{m})>0$ provided $2^{-m}<\mu (K).$ Now choose
compact $C_{m}\subseteq D_{m},$ with $\mu (D_{m}\backslash C_{m})<2^{-m}.$
So $\mu (K\backslash C_{m})<2^{1-m}.$ Also $C_{m}\cap t_{n}C_{m}=\emptyset ,$
for each $n\geq m,$ as $C_{m}\subseteq K;$ but $t_{n}\rightarrow 1_{G},$ so
the compact set $C_{m}C_{m}^{-1}$ contains no interior points. Hence, by
Baire's theorem, neither does $CC^{-1}$ for $C=\bigcup\nolimits_{m}C_{m}$
which differs from $K$ by a null set. $\square $

\bigskip

The significance of Corollary 1 is that in alternative language it asserts a
\textit{Converse Steinhaus-Weil Theorem:}

\bigskip

\noindent \textbf{Proposition 2.} \textit{A regular Borel measure }$\mu $
\textit{on a topological group }$G$ \textit{has the Steinhaus-Weil property
iff either of the following holds:}\newline
\noindent (i) \textit{for each non-null compact subset }$K$\textit{\ the map
}$m_{K}:t\rightarrow \mu (tK)$\textit{\ is subcontinuous at }$1_{G};$\newline
\noindent (ii) \textit{for each non-null compact subset }$K$\textit{\ there
is no `null' sequence }$t_{n}\rightarrow 1_{G}$\textit{\ with }$\mu
(t_{n}K)\rightarrow 0.$

\bigskip

\noindent \textbf{Remark.} In the locally compact case, Simmons and Mospan
both prove that this is equivalent to $\mu $ being absolutely continuous
w.r.t. Haar measure $\eta ;$ see \S 6.2 below. For the more general context
of a Polish group see [BinO9].

\bigskip

For the Baire-category version, which goes back to Piccard and Pettis (see
[BinO6] for references), we recall that the \textit{quasi-interior}, here
conveniently denoted $\tilde{A}$ (or $A^{q}$) of a set $A$ with the Baire
property, is the largest (usual) open set $U$ such that $U\backslash A$ is
meagre; it is a regularly-open set (see [Dug, Ch. 3 Problems, Section 4
Q22], \S 6.3). We note that $(aA)^{q}=aA^{q}.$ We learn from Theorem 4 below
that the counterpart of this `quasi-interior' for measurable sets is
provided by the open sets of the density topology.

\bigskip

\noindent \textbf{Proposition 1B }(\textbf{Baire-Kemperman property) }%
([BinO4, Th. 5.5B/M], [BinO3, \S 5, Th. K]). \textit{In a normed topological
group }$G,$ \textit{Baire under its norm topology, if }$A$\textit{\ is Baire
non-meagre, then there is }$\delta >0$\textit{\ with}%
\[
tA\cap A\in \mathcal{B}_{+}\qquad (||t||<\delta );
\]%
\textit{so in particular,}%
\[
B_{\delta }\subseteq AA^{-1}.
\]

\bigskip

\noindent \textbf{Proof.} If $a\in A^{q},$ the quasi-interior of $A,$ then $%
1_{G}\in a^{-1}\tilde{A},$ which is open. So w.l.o.g. we may take $a=1_{G}.$
Choose $\delta >0$ so that $B:=B_{\delta }\subseteq \tilde{A}.$ Then for $%
t\in B,$ since $t\in tB\cap B,$%
\[
t\tilde{A}\cap \tilde{A}\supseteq tB\cap B\neq \emptyset ,
\]%
so being open, $tB\cap B$ is non-meagre (as $G$ is Baire). But, modulo
meagre sets, $A$ and $\tilde{A}$ are identical. For the remaining assertion,
argue as in the measure case. $\square $

\bigskip

Propositions 1M and 1B are both included in Theorem 4 below in the topology
refinement context with $\mathcal{B}_{+}(\mathcal{D}_{\mathcal{L}})=\mathcal{%
M}_{+}$, and $\mathcal{B}_{+}(\mathcal{D}_{\mathcal{B}})=\mathcal{B}_{+}(%
\mathcal{T}_{d});$ we apply the result in Theorem 5. (Recall the discussion
in \S 2 of their properties, enlisted below.)

\bigskip

\noindent \textbf{Theorem 4 (Displacements Theorem).} \textit{In a
topological group under the topology of }$d_{R}^{G}$\textit{, let }$\mathcal{%
D}$\textit{\ be a topology refining }$d_{R}^{G}$ \textit{and having the
following properties:}

\noindent (i) $\mathcal{D}$\textit{\ is shift-invariant: }$x\mathcal{D}=%
\mathcal{D}$ \textit{for all }$x$\textit{;}

\noindent (ii) $\mathcal{B}_{0}(\mathcal{D})$\textit{\ is left invariant for
all }$x$\textit{;}

\noindent (iii) \textit{`Localization property': }$H\backslash $\textrm{int}$%
_{\mathcal{D}}(H)\in \mathcal{B}_{0}(\mathcal{D})$ \textit{for} $H\in
\mathcal{B}(\mathcal{D})$;

\noindent (iv) \textit{the (left) `Kemperman property': for} $1_{G}\in U\in
\mathcal{D}$ \textit{there is }$\delta =\delta _{U}>0$\textit{\ with}%
\[
tU\cap U\in \mathcal{B}_{+}(\mathcal{D})\qquad (||t||<\delta );
\]

\noindent (v) $\mathcal{D}$ \textit{is a strong generalized Gandy-Harrington
topology: modulo} $\mathcal{B}_{0}(\mathcal{D})$ \textit{each} $\mathcal{B}(%
\mathcal{D})$ \textit{set is analytic under} $d_{R}^{X},$ \textit{so that }$%
\mathcal{D}$\textit{\ is a Baire space.}

\textit{Then for }$A\in \mathcal{B}_{+}(\mathcal{D})$ \textit{and quasi all }%
$a\in A,$\textit{\ there is }$\varepsilon =\varepsilon _{A}(a)>0$ \textit{%
with}%
\[
axa^{-1}A\cap A\in \mathcal{B}_{+}(\mathcal{D})\qquad (||x||<\varepsilon );
\]%
\textit{so in particular, with }$\gamma _{a}(x):=axa^{-1},$%
\[
\gamma _{a}(B_{\varepsilon })\subseteq AA^{-1}\qquad \text{off a }\mathcal{B}%
_{0}(\mathcal{D})\text{-set of }a\in A.
\]%
\

\noindent \textbf{Proof.} For some $U\in \mathcal{D}$ and $N,A^{\prime }\in
\mathcal{B}_{0}(\mathcal{D}),$ $A\backslash A^{\prime }=U\backslash N;$ by
(i) and (ii), $U_{a}=a^{-1}U\in \mathcal{D}$ and $N_{a}=a^{-1}N\in \mathcal{B%
}_{0}(\mathcal{D}),$ so $a^{-1}A\supseteq U_{a}\backslash N_{a}$ for $a\in
A\backslash A^{\prime }.$ By the Kemperman property, for $a\in A\backslash
A^{\prime }$ there is $\varepsilon =\varepsilon _{A}(a)>0$ with
\[
xU_{a}\cap U_{a}\in \mathcal{B}_{+}(\mathcal{D})\qquad (||x||<\varepsilon ).
\]%
Fix $a\in A\backslash A^{\prime }.$ Working modulo $\mathcal{B}_{0}(\mathcal{%
D})$-sets which are left invariant (by the Nikodym property), $%
xa^{-1}A\supseteq xU_{a}$ and%
\[
a^{-1}A\cap xa^{-1}A\supseteq U_{a}\cap (xU_{a}),\qquad (||x||<\varepsilon
).
\]%
So, since also $axa^{-1}A\in \mathcal{B}(\mathcal{D}),$ again by the Nikodym
property%
\[
axa^{-1}A\cap A\in \mathcal{B}_{+}(\mathcal{D})\qquad (||x||<\varepsilon ).
\]%
For the remaining assertion, argue as in Prop. 1M above. $\square $

\bigskip

\noindent \textbf{Remark.} The map $\gamma _{a}:x\mapsto axa^{-1}$ in the
preceding theorem is a homeomorphism under the topology of the topological
group $d_{R}^{G}$ (being continuous, with continuous inverse $\gamma
_{a^{-1}})$, so $\gamma _{a}(B_{\varepsilon _{A}(a)})$ is open, and also
open under the finer topology $\mathcal{D}$.

\bigskip

Theorem 4 above identifies additional topological properties enabling the
Kemperman property to imply the Steinhaus-Weil property. So we regard it as
a \textit{weak} Steinhaus-Weil property; the extent to which it is weaker is
clarified by Proposition 3 below. For this we need a $\sigma $-algebra.
Recall that $E\subseteq G$ is \textit{universally measurable} ($E\in
\mathcal{U}(G)$) if $E$ is measurable with respect to every measure $\mu \in
\mathcal{P}(G)$ -- for background, see e.g. [Kec, \S 21D], cf. [Fre, 434D,
432]; these form a $\sigma $\textit{-algebra}. Examples are analytic subsets
(see e.g. [Rog, Part 1 \S 2.9], or [Kec, Th. 21.10], [Fre, 434Dc]) and the $%
\sigma $-algebra that they generate. Beyond these are the \textit{provably }$%
\mathbf{\Delta }_{2}^{1}$ sets of [FenN].

\bigskip

\noindent \textbf{Proposition 3 }(cf. BinO8, Cor. 2])\textbf{.} \textit{For}
$\mathcal{H}$ \textit{a left invariant }$\sigma $\textit{-ideal in }$G$%
\textit{, put}%
\[
\mathcal{U}_{+}(\mathcal{H})=\mathcal{U}(G)\backslash \mathcal{H}.
\]%
\textit{If }$\mathcal{U}_{+}(\mathcal{H})$ \textit{has the Steinhaus-Weil
property, then }$\mathcal{U}_{+}(\mathcal{H})$ \textit{has the Kemperman
property.}

\bigskip

\noindent \textbf{Proof.} Suppose otherwise. Then, for some $E\in \mathcal{U}%
_{+}(\mathcal{H})$ and some $t_{n}\rightarrow 1_{G},$ each of the sets $%
E\cap t_{n}E$ is in $\mathcal{H}$, and so%
\[
E_{0}:=E\backslash \dbigcup\nolimits_{n}t_{n}E=E\backslash
\dbigcup\nolimits_{n}(E\cap t_{n}E)\in \mathcal{U}(G),
\]%
which is in $\mathcal{U}_{+}(\mathcal{H})$ (as otherwise $E=E_{0}\cup $ $%
\dbigcup [E\cap t_{n}E]\in \mathcal{H}$). So $1_{G}\in \mathrm{int}%
(E_{0}E_{0}^{-1}).$ So, for ultimately all $n,$ $t_{n}\in E_{0}E_{0}^{-1},$
and then $E_{0}\cap t_{n}E_{0}\neq \emptyset .$ But, as $E_{0}\subseteq E,$ $%
E_{0}\cap t_{n}E_{0}\subseteq E_{0}\cap t_{n}E=\emptyset ,$ a contradiction.
$\square $

\bigskip

Of interest above is the case $\mathcal{H}$ of \textit{left-Haar-null} sets
[Sol2] of a Polish group $G$ (cf. \S 6.5,6). We close with a result on the
density topology $\mathcal{D}_{\mu }$ generated in a Polish group $G$ from
an \textit{atomless} measure $\mu $; this is an immediate corollary of
[Mar]. Unlike the Hashimoto ideal topologies, which need not be Baire
topologies, these are Baire. The question of which sets have the
Steinhaus-Weil property under $\mu ,$ hinges on the choice of $\mu $ (see
the earlier cautionary example of this section), and indeed on further
delicate subcontinuity considerations, related to [Sol2], for which see
[BinO9]. In this connection see [Oxt1] and [DieS, Ch. 10].

\bigskip

\noindent \textbf{Proposition 4.} \textit{For }$G$\textit{\ a Polish group
with metric topology }$\mathcal{T}_{d}$,\textit{\ }$\mathcal{\beta }$
\textit{a countable basis, and atomless} $\mu \in \mathcal{P}(G)$ \textit{%
with }$\mu (U)>0$ \textit{for all non-empty }$U\in \beta $, \textit{then,
with density at }$g\in G$\textit{\ computed by reference to }$\beta _{g}$ $%
:=\{B\in \beta :g\in B\}$\textit{, the Lebesgue density theorem holds for }$%
\mu .$ \textit{So the generated density topology }$\mathcal{D}_{\mu }$
\textit{refines }$\mathcal{T}_{d}$ \textit{and is a Baire topology.}

\bigskip

\noindent \textbf{Proof.} As $\mathcal{\beta }$ comprises open sets, the
Vitali covering lemma applies (see [Bru, \S 6.3-4]), and implies that the
measure $\mu $ obeys a density theorem (that $\mu $-almost all points of a
measurable set are density points -- cf. [Kuc, Th. 3.5.1]). So by results of
Martin [Mar, Cor. 4.4], the family $\mathcal{D}_{\mu }$ of sets all of whose
points have (outer) density 1 are $\mu $-measurable, and form a topology on $%
G$ [Mar, Th. 4.1]$.$ It is a Baire topology, by [Mar, Cor. 4.13]. It refines
$\mathcal{T}_{d}$: the points of any open set have density 1, because the
differentiation basis consists of open sets, and these all are $\mu $%
-non-null. $\square $

\bigskip

\noindent \textbf{Remark.} The assumption of regularity subsumed in $\mu \in
\mathcal{P}(G)$ is critical; in its absence the density theorem may fail:
see [Kha, Ch.8 Th. 1], where for $\mu $ non-regular there is a $\mu $%
-measurable set with just one density point.

\section{A Shift Theorem}

Theorem 5 below establishes a compactness-like property of a density
topology $\mathcal{D}$ which, according to Corollary 2 below, implies the
Steinhaus-Weil property for sets in $A\in \mathcal{B}_{+}(\mathcal{D}).$ So
we may regard it as a strong Steinhaus-Weil-like property. Its prototype
arises in the relevant infinite combinatorics (the Kestelman-Borwein-Ditor
Theorem, KBD: see [BinO6], [Ost3]). The setting for the theorem is that of
the Displacements Theorem (Th. 4 above), a key ingredient of which is the
Kemperman property, a weak Steinhaus-Weil-like property. So Theorem 5
establishes the equivalence of a strong and a weak Steinhaus-Weil-like
property in the presence of additional topological restrictions on the
relevant refinement topology: invariance under shift, localization and some
analyticity (namely, a weak base of analytic sets). We raise and leave open
the question as to whether the three topological restrictions listed above
are minimal here.

\bigskip

\noindent \textbf{Theorem 5 (Fine Topology Shift Theorem).} \textit{In a
topological group under the topology of }$d_{R}^{G}$\textit{, let }$\mathcal{%
D}$\textit{\ be a topology refining }$d_{R}^{G}$ \textit{and having the
following properties:}

\noindent (i) $\mathcal{D}$\textit{\ is shift-invariant: }$x\mathcal{D=D}$
\textit{for all }$x$\textit{;}

\noindent (ii) \textit{`Localization property': }$H\backslash $\textrm{int}$%
_{\mathcal{D}}(H)\in \mathcal{B}_{0}(\mathcal{D})$ \textit{for} $H\in
\mathcal{B}(\mathcal{D})$;

\noindent (iii) \textit{the (left) `Kemperman property': for} $1_{G}\in U\in
\mathcal{D}$ \textit{there is }$\delta =\delta _{U}>0$\textit{\ with}%
\[
U\cap (tU)\in \mathcal{B}_{+}(\mathcal{D})\qquad (||t||<\delta );
\]

\noindent (iv) $\mathcal{D}$ \textit{is a strong generalized
Gandy-Harrington topology: modulo} $\mathcal{B}_{0}(\mathcal{D})$ \textit{%
each} $\mathcal{B}(\mathcal{D})$ \textit{set is analytic under} $d_{R}^{G},$
\textit{so that }$\mathcal{D}$\textit{\ is a Baire space.}

\textit{Then, for }$z_{n}\rightarrow 1_{G}$\textit{\ (\textquotedblleft null
sequence\textquotedblright ) and }$A\in \mathcal{B}_{+}(\mathcal{D}),$%
\textit{\ for quasi all }$a\in A$\textit{\ there is an infinite set} $%
\mathbb{M}_{a}$ \textit{such that }%
\[
\{z_{m}a:m\in \mathbb{M}_{a}\}\subseteq A.
\]%
\noindent \textbf{Proof of Theorem 5.} Since the asserted property is
monotonic, we may assume by the strong Gandy-Harrington property that $A$ is
analytic in the topology of $d_{R}^{G}$. So write $A=K(I)$ with $K\ $%
upper-semicontinuous and single-valued. Below, for greater clarity we write $%
B(x,r)$ for the open $r$-ball centered at $x.$ For each $n\in \omega $ we
find inductively integers $i_{n},$ points $x_{n},y_{n},a_{n}$ with $a_{n}\in
A$, numbers $r_{n}\downarrow 0,s_{n}\downarrow 0,$ analytic subsets $A_{n}$
of $A,$ and closed nowhere dense sets $\{F_{m}^{n}:m\in \omega \}$ w.r.t. $%
\mathcal{D}$ and $D_{n}\in \mathcal{D}$ such that:
\[
K(i_{1},...,i_{n})\supseteq a_{n}x_{n}a_{n}^{-1}A_{n}\cap A_{n}\in \mathcal{B%
}_{+}(\mathcal{D}),
\]%
\[
K(i_{1},...,i_{n})\supseteq a_{n}x_{n}a_{n}^{-1}A_{n}\cap A_{n}\supseteq
D_{n}\cap B(y_{n},s_{n})\backslash \bigcup\nolimits_{m\in \omega }F_{m}^{n},
\]%
\[
y_{n}\in B(a_{n}x_{n},r_{n})\text{ and }D_{n}\cap B(y_{n},s_{n})\cap
\bigcup\nolimits_{m,k<n}F_{m}^{k}=\emptyset .
\]%
Assuming this done for $n,$ since $K(i_{1},...,i_{n})=\bigcup%
\nolimits_{k}K(i_{1},...,i_{n},k)$ is non-negligible, there is $i_{n+1}$
with $K(i_{1},...,i_{n},i_{n+1})\cap a_{n}x_{n}a_{n}^{-1}A_{n}\cap A_{n}$
non-negligible. Put $A_{n+1}:=K(i_{1},...,i_{n},i_{n+1})\cap
a_{n}x_{n}a_{n}^{-1}A_{n}\cap A_{n}\subseteq K(i_{1},...,i_{n+1}).$ As $%
A_{n+1}$ is non-negligible, we may pick $a_{n+1}\in A_{n+1}$ and $%
\varepsilon (a_{n+1},A_{n+1})$ as in Theorem 4 above; also pick $m(n)$ so
large that $||z_{m}||<\varepsilon (a_{n+1},A_{n+1})$ for $m\geq m(n)$ and
that for $x_{n+1}:=a_{n+1}z_{m(n)}^{-1}a_{n+1}^{-1}$ also $||x_{n+1}||\leq
2^{-n}$ (the latter is possible, as $G$ is a topological group in the
group-norm topology). Then
\[
K(i_{1},...,i_{n+1})\supseteq a_{n+1}x_{n+1}a_{n+1}^{-1}A_{n+1}\cap
A_{n+1}\in \mathcal{B}_{+}(\mathcal{D}).
\]%
So by the Banach Category Theorem, for some $D_{n+1}\in \mathcal{D}$ and
some positive $r_{n+1}<r_{n}/2$%
\[
K(i_{1},...,i_{n+1})\supseteq a_{n+1}x_{n+1}a_{n+1}^{-1}A_{n+1}\cap
A_{n+1}\supseteq D_{n+1}\cap B(a_{n+1}x_{n+1},r_{n+1})\backslash
\bigcup\nolimits_{m\in \omega }F_{m}^{n+1},
\]%
for some closed nowhere dense sets $\{F_{m}^{n+1}:m\in \omega \}$ of $%
\mathcal{D}$. Since the set $\bigcup\nolimits_{m,k<n+1}F_{m}^{k}$ is closed
and nowhere dense, there is $y_{n+1}\in D_{n+1}\cap
B(a_{n+1}x_{n+1},r_{n+1}) $ and positive $s_{n+1}<s_{n}/2$ so small that $%
B(y_{n+1},s_{n+1})\subseteq B(a_{n+1}x_{n+1},r_{n+1})$ and $D_{n+1}\cap
B(y_{n+1},s_{n+1})\cap \bigcup\nolimits_{m,k<n+1}F_{m}^{k}=\emptyset .$ So%
\[
D_{n+1}\cap B(a_{n+1}x_{n+1},r_{n+1})\backslash \bigcup\nolimits_{m\in
\omega }F_{m}^{n+1}\supseteq D_{n+1}\cap B(y_{n+1},s_{n+1})\backslash
\bigcup\nolimits_{m\in \omega }F_{m}^{n+1},
\]%
completing the induction.

By the Analytic Cantor Theorem [Ost1, Th. AC, Section 2], there is $t$ with
\[
\{t\}=K(i)\cap \bigcap\nolimits_{n}B(y_{n},s_{n})\subseteq
\bigcap\nolimits_{n}a_{n}x_{n}a_{n}^{-1}A_{n}\cap A_{n}.
\]%
So $t\in A.$ Fix $n;$ then $t\notin \bigcup\nolimits_{m\in \omega }F_{m}^{n}$
(since $D_{m+1}\cap B(y_{m+1},s_{m+1})\cap
\bigcup\nolimits_{k<m}F_{k}^{n}=\emptyset $ for each $m),$ and so $t\in T$
and%
\[
t\in D_{n}\cap B(y_{n},s_{n})\backslash \bigcup\nolimits_{m\in \omega
}F_{m}^{n}\subseteq a_{n}x_{n}a_{n}^{-1}A_{n}\cap A_{n}\subseteq
K(i_{1},...,i_{n})\subseteq A.
\]%
As $t\in a_{n+1}x_{n+1}a_{n+1}^{-1}A_{n+1},$ $%
a_{n+1}x_{n+1}^{-1}a_{n+1}^{-1}t=z_{m(n)}t\in A_{n+1}\subseteq A.$ So $%
\{z_{m(n)}t:n\in \omega \}\subseteq A.$

Now $d_{R}^{G}(a_{n}x_{n},t)\leq
d_{R}^{G}(a_{n}x_{n},y_{n})+d_{R}^{G}(y_{n},t)\rightarrow 0$, so $%
a_{n}\rightarrow t,$ since $x_{n}\rightarrow 1_{G}$ and so%
\[
d_{R}^{G}(a_{n},t)=d_{R}^{G}(a_{n}x_{n},tx_{n})\leq
d_{R}^{G}(a_{n}x_{n},t)+d_{R}^{G}(1,tx_{n}t^{-1})\rightarrow 0,
\]%
again as $G$ is a topological group under the $d_{R}^{G}$-topology. $\square
$

\bigskip

\noindent \textbf{Corollary 2.} \textit{In the setting of Theorem 5 above,
the sets of }$\mathcal{B}_{+}(\mathcal{D})$\textit{\ have the Steinhaus-Weil
property (\S 1).}

\bigskip

\noindent \textbf{Proof }(cf. [Sol1, Th. 1(ii)], [BinO4, Th. 6.5]).
Otherwise, for some set $A\in \mathcal{B}_{+}(\mathcal{D})$ we may select $%
z_{n}\notin AA^{-1}$ with $||z_{n}||<1/n.$ Then there are $a\in A$ and a
subsequence $m(n)$ with $z_{m(n)}a\in A;$ so $z_{m(n)}\in AA^{-1},$ a
contradiction. $\square $

\section{Complements}

\noindent 1. \textit{Weil topology }([Wei, Ch. VII, \S 31], cf. [Hal, Ch.
XII \S 62]). We recall that for $G$ a group with a $\sigma $-finite
left-invariant measure $|.|$ on a $\sigma $-ring $\mathcal{M}$ of
left-invariant sets and $(x,y)\mapsto (x,xy)$ measurability-preserving, the
\textit{Weil topology }is generated by the family of pseudo-norms%
\[
||g||_{E}:=|gE\triangle E|,
\]%
for $E\in \mathcal{M}_{+}$ (with $\mathcal{M}_{+}$ the family of measurable
sets with finite positive measure), so that $||g||_{E}\leq |E|$. Provided
the pseudo-norms are \textit{separating} (i.e. $||g||_{E}>0$ for any $g\neq
1_{G}$ and some $E\in \mathcal{M}_{+}$ as in (iii) above)$,$ $G$ is a
topological group under the Weil topology [Hal, 62E]; equivalently, the
topology is generated by the neighbourhood base $\mathcal{N}%
_{1}:=\{DD^{-1}:1_{G}\in D\in \mathcal{M}^{+}\}$, reminiscent of the
Steinhaus-Weil Theorem. The proof relies on a kind of fragmentation lemma
(see [BinO9]). That in turn depends on Fubini's Theorem [Hal, 36C], via the
average theorem [Hal, 59.F]:%
\[
\int_{G}|g^{-1}A\cap B|dg=|A|\cdot |B^{-1}|\qquad (A,B\in \mathcal{M}),
\]%
($g=ab^{-1}$ iff $g^{-1}a=b),$ and may be interpreted as demonstrating the
continuity at $1_{G}$ of $||.||_{E}$ under the density topology.

\noindent 2. \textit{Steinhaus-Weil property of a Borel measure. }In a
locally compact group $G,$ the family $\mathcal{M}_{+}(\mu )$ of finite
non-null measurable sets of a Borel measure $\mu $ on $G$ fails to have the
Steinhaus-Weil property iff there are a null sequence $z_{n}\rightarrow
1_{G} $ and a non-null compact set $K$ with $\lim_{n}\mu (t_{n}K)=0,$ as
observed by Mospan [Mos] (in $\mathbb{R}).$ Equivalently, this is so iff the
measure $\mu $ is not absolutely continuous with respect to Haar measure --
cf. [Sim] and [BinO9].\newline
\noindent 3. \textit{Regular open sets. }Recall that $U\ $is regular open if
$U=$\textrm{int}$($\textrm{cl}$U),$ and that \textrm{int}$($\textrm{cl}$U)$
is itself regular open; for background see e.g. [GivH, Ch. 10], or [Dug, Ch.
3 Problems, Section 4 Q22]. For $\mathcal{D=D}_{\mathcal{B}}$ the
Baire-density topology of a normed topological group, let $\mathcal{D}_{%
\mathcal{B}}^{RO}$ denote the regular open sets. For $D\in \mathcal{D}_{%
\mathcal{B}}^{RO}$, put%
\[
N_{D}:=\{t\in G:tD\cap D\neq \emptyset \}=DD^{-1},\qquad \mathcal{N}%
_{1}:=\{N_{D}:1_{G}\in D\in \mathcal{D}_{RO}\};
\]%
then $\mathcal{N}_{1}$ is a base at $1_{G}$ (since $1_{G}\in C\in \mathcal{D}%
_{RO}$ and $1_{G}\in D\in \mathcal{D}_{RO}$ yield $1_{G}\in C\cap D\in
\mathcal{D}_{RO})$ comprising $\mathcal{T}$-neighbourhoods that are $%
\mathcal{D}_{\mathcal{B}}$-open (since $DD^{-1}=\bigcup \{Dd^{-1}:d\in D\}$.
We raise the (metrizability)\ question, by analogy with the Weil topology of
a measurable group: with $\mathcal{D}_{\mathcal{B}}$ above replaced by a
general density topology $\mathcal{D}$ on a group $G,$ when is the topology
generated by $\mathcal{N}_{1}$ on $G$ a norm topology? Some indications of
an answer may be found in [ArhT, \S 3.3]. We note the following plausible
answer: if there exists a separating sequence $D_{n},$ i.e. such that for
each $g\neq 1_{G}$ there is $n$ with $||g||_{D_{n}}=1,$ then%
\[
||g||:=\sum\nolimits_{n}2^{-n}||g||_{D_{n}}
\]%
is a norm, since it is separating and, by the Nikodym property, $(D\cap
g^{-1}D)=g^{-1}(gD\cap D)\in \mathcal{B}_{0}$.

\noindent 4. \textit{The Effros Theorem }asserts that a transitive
continuous action of a Polish group $G$ on a space $X$ of second category in
itself is necessarily `open', or more accurately is microtransitive (the
(continuous) evaluation map $e_{x}:g\mapsto g(x)$ takes open nhds $E$ of $%
1_{G}$ to open nhds that are the orbit sets $E(x)$ of $x$). It emerges that
this assertion is very close to the shift-compactness property: see [Ost6].
The Effros Theorem reduces to the Open Mapping Theorem when $G,X\ $are
Banach spaces regarded as additive groups, and $G$ acts on $X\ $by way of a
linear surjection $L:G\rightarrow X$ via $g(x)=L(g)+x.$ Indeed, here $%
e_{0}(E)=L(E).$ For a neat proof, choose an open neighbourhood $U$ of $0$ in
$G$ with $E\supseteq U-U$; then $L(U)$ is Baire (being analytic) and
non-meagre (since $\{L(nU):n\in N\}$ covers $X),$ and so $L(U)-L(U)\subseteq
L(E)$ is an open nhd of $0$ in $X.$

\noindent 5. \textit{Amenability at 1}. Solecki defines $G$ to be amenable
at 1 if given $\mu _{n}\in \mathcal{P}(G)$ with $1_{G}\in $ \textrm{supp}$%
(\mu _{n})$ there are $\nu $ and $\nu _{n}$ in $\mathcal{P}(G)$ with
\[
\nu _{n}\ast \nu (K)\rightarrow \nu (K)\qquad (K\in \mathcal{K}(G)).
\]%
(The origin of the term may be traced to a localization, via the restriction
of supports to contain $1_{G},$ of a \textit{Reiter-like condition }[Pat,
Prop. 0.4] characterizing amenability itself.) It is proved in [Sol2, Th.
1(ii)] that, in the class of Polish groups $G$ that are amenable at $1_{G},$
the Steinhaus-Weil property holds for universally measurable sets that are
not left-Haar-null$;$ this includes Polish abelian groups [Sol2, Prop. 3.3].
The relativized notion of subcontinuity: \textit{on a compact }$K$ \textit{%
along a null sequence} $\{t_{n}\}$ (which requires some subsequence $\mu
(t_{m(n)}K)$ to be bounded away from $0,$ provided $\mu (K)>0$) yields a
connection to amenability at $1_{G}$, which we explore elsewhere [BinO9].

\noindent 6. \textit{Haar-null and left-Haar-null. }The two families
coincide in Polish abelian groups, and in locally compact second countable
groups (where they also coincide with the sets of Haar measure zero -- by an
application of the Fubini theorem). The former family, however, is in
general smaller; indeed, (universally measurable) non-Haar-null sets need
not have the Steinhaus-Weil property, whereas the (universally measurable)
non-left-Haar-null sets do -- see [Sol2].

\noindent 7. \textit{Left-Haar-null Kemperman property. }We note, as this is
thematic, that the family of (universally measurable) non-left-Haar-null
sets has the left Kemperman property [BinO9, Lemma 1].

\noindent 8. \textit{Beyond local compactness:\ Haar category-measure
duality. }In the absence of Haar measure, the definition of left-Haar-null
subsets of a topological group $G$ requires $\mathcal{U}(G),$ the
universally measurable sets -- by dint of the role of the totality of
(probability) measures on $G$. The natural dual of $\mathcal{U}(G)$ is the
class $\mathcal{U}_{\mathcal{B}}(G)$ of \textit{universally Baire sets},
defined,for $G$ with a Baire topology, as those sets $B$ whose preimages $%
f^{-1}(B)$ are Baire in any compact Hausdorff space $K$ for any continuous $%
f:K\rightarrow G$. Initially considered in [FenMW] for $G=\mathbb{R}$, these
have attracted continued attention for their role in the investigation of
axioms of determinacy and large cardinals -- see especially [Woo]; cf.
[MarS].

Analogously to the left-Haar-null sets, define in $G$ the family of \textit{%
left-Haar-meagre} sets, $\mathcal{HM}(G)$, to comprise the sets $M$
coverable by a universally Baire set $B$ for which there are a compact
Hausdorff space $K$ and a continuous $f:K\rightarrow G$ with $f^{-1}(gB)$
meagre in $K$ for all $g\in G.$ These were introduced, in the \textit{abelian%
} Polish group setting with $K$ metrizable, by Darji [Dar], cf. [Jab1], and
shown there to form a $\sigma $-ideal of meagre sets (co-extensive with the
meagre sets for $G\ $locally compact); as $\mathcal{HM}(G)\mathcal{\subseteq
B}_{0}(G),\mathcal{\ }$the family is not studied here.

\noindent 9. \textit{Steinhaus }$AA^{-1}$\textit{\ and }$AB^{-1}$\textit{\
properties.\protect\footnote{ See the Appendix of this arXiv version for a fuller account.}} If the subsets of $G$ lying in a family $\mathcal{H}$ have the
property that $AA^{-1}$ for $A\in \mathcal{H}$ has non-empty $\tau $%
-interior, for $\tau $ a translation invariant topology, and furthermore, as
in the Haar-Kemperman property, for $A,B\in \mathcal{H}$ there is $g\in G$
such that $C:=gA\cap B\in \mathcal{H},$ then of course $g^{-1}CC^{-1}%
\subseteq AB^{-1},$ and so the latter has non-empty $\tau $-interior. By the
Average Theorem (\S 6.1 above), this is the case for $G$ locally compact
with $\tau =\mathcal{T}_{d}$ and $\mathcal{H=L}_{+}$ the Haar-measurable
non-null sets [Hal, \S 59F] (cf. [TomW, \S 11.3], and [BinO5] for $G=\mathbb{%
R}$); other examples of families $\mathcal{H}$ are provided by certain
refinement topologies $\tau $ -- see [BinO9]. However, Mato\u{u}skov\'{a}
and Zelen\'{y} [MatZ] show that in any non-locally compact abelian Polish
group there are closed non (left) Haar null sets $A,B$ such that $A+B$ has
empty interior. Recently, Jab\l o\'{n}ska [Jab2] has shown that likewise in
any non-locally compact abelian Polish group there are closed non-Haar
meager sets $A,B$ such that $A+B$ has empty interior.

\noindent 10. \textit{Non-separability. }The links between the Effros
theorem above, the Baire theorem and the Steinhaus-Weil theorem are pursued
at length in [Ost6]. There, any separability assumption is avoided. Instead
\textit{sequential} methods are used, for example shift-compactness
arguments.

\noindent 11. \textit{Metrizability and Christensen's Theorem.} In
connection with the role of analyticity in the generalized Gandy-Harrington
property of \S 2, note that an analytic topological group is metrizable; so
if it is also a Baire space, then it is a Polish group [HofT, Th. 2.3.6].

\noindent 12. \textit{Strong Kemperman property: qualitative versus
quantitative measure theory. }We note that property (iv) of Theorem 1
corresponds to the following quantitative, linear Lebesgue-measure property,
which we may name the \textit{strong Kemperman property} (see [Kem], [Kuc,
Lemma 3.7.2]):

\noindent (iv)* for $0\in U\in \mathcal{D}_{\mathcal{L}}$ there is $\delta
>0 $ such that for all $|t|<\delta $
\[
|U\cap (t+U)|\geq \varepsilon .
\]%
This is connected with the continuity of a Weil-like group norm on $(\mathbb{%
R},+)$. Indeed, since%
\[
|U\cap (t+U)|=|U|-|U\triangle (t+U)|/2,
\]%
the inequality above is equivalent to
\[
||t||_{U}:=|U\triangle (t+U)|\leq 2(|U|-\varepsilon ).
\]%
The latter holds for any $0<\varepsilon <|U|\ $and for sufficiently small $%
t, $ by the continuity of the norm $||t||_{U}$.

\bigskip

\textbf{Acknowledgements. }We thank the Referee for a careful and scholarly
reading of the paper, and for some very useful presentational
suggestions.\newpage

\textbf{References.}

\noindent \lbrack ArhT] A. Arhangelskii, M. Tkachenko, \textsl{Topological
groups and related structures}. Atlantis Studies in Mathematics, 1. Atlantis
Press, Paris; World Scientific Publishing Co. Pte. Ltd., Hackensack, NJ,
2008. \newline
\noindent \lbrack BalRS] M. Balcerzak, A. Ros\l anowski, S. Shelah, Ideals
without ccc, J\textsl{. Symbolic Logic}, \textbf{63} (1998), 128-148.\newline
\noindent \lbrack BarFN] A. Bartoszewicz, M. Filipczak, T. Natkaniec, On
Smital properties. \textsl{Topology Appl.} \textbf{158} (2011), 2066--2075.%
\newline
\noindent \lbrack BeeV] G. Beer, L. Villar, Borel measures and Hausdorff
distance. \textsl{Trans. Amer. Math. Soc.} \textbf{307} (1988), 763--772.%
\newline
\noindent \lbrack BinGT] N. H. Bingham, C. M. Goldie and J. L. Teugels,
\textsl{Regular variation}, 2nd ed., Cambridge University Press, 1989 (1st
ed. 1987).\newline
\noindent \lbrack BinO1] N. H. Bingham and A. J. Ostaszewski, Beyond
Lebesgue and Baire: generic regular variation, \textsl{Coll. Math. }\textbf{%
116} (2009), 119-138.\newline
\noindent \lbrack BinO2] N. H. Bingham and A. J. Ostaszewski, Beyond
Lebesgue and Baire II: Bitopology and measure-category duality, \textsl{%
Coll. Math.}, \textbf{121} (2010), 225-238.\newline
\noindent \lbrack BinO3] N. H. Bingham and A. J. Ostaszewski, Kingman,
category and combinatorics. \textsl{Probability and Mathematical Genetics}
(Sir John Kingman Festschrift, ed. N. H. Bingham and C. M. Goldie), 135-168,
London Math. Soc. Lecture Notes in Mathematics \textbf{378}, CUP, 2010.
\newline
\noindent \lbrack BinO4] N. H. Bingham and A. J. Ostaszewski, Normed groups:
Dichotomy and duality. \textsl{Dissert. Math.} \textbf{472} (2010), 138p.
\newline
\noindent \lbrack BinO5] N. H. Bingham and A. J. Ostaszewski, Regular
variation without limits, \textsl{J. Math. Anal. Appl.} \textbf{370} (2010),
322-338.\newline
\noindent \lbrack BinO6] N. H. Bingham and A. J. Ostaszewski, Dichotomy and
infinite combinatorics: the theorems of Steinhaus and Ostrowski. \textsl{%
Math. Proc. Camb. Phil. Soc.} \textbf{150} (2011), 1-22. \newline
\noindent \lbrack BinO7] N. H. Bingham and A. J. Ostaszewski, The Steinhaus
theorem and regular variation: de Bruijn and after, \textsl{Indag. Math.}
\textbf{24} (2013), 679-692.\newline
\noindent \lbrack BinO8] N. H. Bingham and A. J. Ostaszewski,
Category-measure duality: convexity, mid-point convexity and Berz
sublinearity, \textsl{Aequationes Math.}, \textbf{91.5} (2017), 801--836 (
fuller version: arXiv1607.05750).\newline
\noindent \lbrack BinO9] N. H. Bingham and A. J. Ostaszewski, The
Steinhaus-Weil property: its converse, Solecki amenability and
subcontinuity, arXiv1607.00049v3.\newline
\noindent \lbrack Bog] V.I. Bogachev, \textsl{Gaussian Measures}, Math.
Surveys \& Monographs \textbf{62}, Amer Math Soc., 1998.\newline
\noindent \lbrack Bor] K. C. Border, \textsl{Fixed point theorems with
applications to economics and game theory.} 2$^{\text{nd}}$ ed. Cambridge
University Press, 1989 (1$^{\text{st}}$ ed.1985).\newline
\noindent \lbrack Bou] A. Bouziad, Continuity of separately continuous group
actions in p-spaces, \textsl{Topology Appl.} \textbf{71}, (1996),119-124%
\newline
\noindent \lbrack Bru] A. M. Bruckner, Differentiation of integrals. \textsl{%
Amer. Math. Monthly} \textbf{78}.9 (1971), 1-51\newline
\noindent \lbrack Chr1] J. P. R. Christensen, On sets of Haar measure zero
in abelian Polish groups. Proceedings of the International Symposium on
Partial Differential Equations and the Geometry of Normed Linear Spaces
(Jerusalem, 1972). \textsl{Israel J. Math.} \textbf{13} (1973), 255--260.%
\newline
\noindent \lbrack Chr2] J. P. R. Christensen, \textsl{Topology and Borel
structure. Descriptive topology and set theory with applications to
functional analysis and measure theory.} North-Holland Mathematics Studies,
Vol. 10, 1974.\newline
\noindent \lbrack CichKW] J. Cicho\'{n}, A. Kharazishvili, B. W\k{e}glorz,
\textsl{Subsets of the real line}, U. \L \'{o}dzki,1995.\newline
\noindent \lbrack CieJ] K. Ciesielski, J. Jasi\'{n}ski, Topologies making a
given ideal nowhere dense or meager. \textsl{Topology Appl.} \textbf{63}
(1995), 277--298.\newline
\noindent \lbrack CieLO] K. Ciesielski, L. Larson, K. Ostaszewski, $\mathcal{%
I}$\textsl{-density continuous functions}, Amer. Math. Soc. Mem. 107, 1994,
133p.\newline
\noindent \lbrack Dan] Dang Anh Tuan, A short proof of the converse to a
theorem of Steinhaus-Weil, arXiv:1511.05657v1.\newline
\noindent \lbrack Dar] U. B. Darji, On Haar meager sets. \textsl{Topology
Appl.} \textbf{160} (2013), 2396--2400.\newline
\noindent \lbrack Den] A. Denjoy, Sur les fonctions deriv\'{e}es sommables,
\textsl{Bull. Soc. Math. Fr.} \textbf{43} (1915), 161-248.\newline
\noindent \lbrack DieS] J. Diestel, A. Spalsbury, \textsl{The joys of Haar
measure.} Grad. Studies in Math. \textbf{150}. Amer. Math. Soc., 2014.%
\newline
\noindent \lbrack Dug] J. Dugundji, \textsl{Topology}, Allyn \& Bacon, 1966.%
\newline
\noindent \lbrack EvaG] L.C. Evans, R.F. Gariepy, \textsl{Measure theory and
fine properties of functions}, CRC, 1992.\newline
\noindent \lbrack FenMW] Q. Feng, M. Magidor, H. Woodin, Universally Baire
sets of reals, in H. Judah, W. Just, H. Woodin (eds.), \textsl{Set theory of
the continuum}, 203--242, Math. Sci. Res. Inst. Publ. \textbf{26}, Springer,
1992.\newline
\noindent \lbrack FenN] J. E. Fenstad, D. Normann, On absolutely measurable
sets. \textsl{Fund. Math.} \textbf{81} (1973/74), no. 2, 91--98.\newline
\noindent \lbrack FilW] M. Filipczak, W. Wilczy\'{n}ski, Strict density
topology on the plane. Measure case. \textsl{Rend. Circ. Mat. Palermo }(2)%
\textbf{\ 60} (2011),113--124.{\normalsize \newline
}\noindent \lbrack Fre] D. Fremlin, \textsl{Measure theory Vol. 3: Measure
algebras.} Corrected 2$^{\text{nd}}$second printing of the 2002 original.
Torres Fremlin, Colchester, 2004.\newline
\noindent \lbrack Ful] F. V. Fuller, Relations among continuous and various
non-continuous functions. \textsl{Pacific J. Math.} \textbf{25} (1968),
495--509.\newline
\noindent \lbrack Gao] Su Gao, \textsl{Invariant descriptive set theory}.
Pure and Applied Mathematics \textbf{293}, CRC Press, 2009.\newline
\noindent \lbrack GivH] S. Givant and P. Halmos, \textsl{Introduction to
Boolean Algebras}, Springer, 2009 (1$^{\text{st}}$ ed., P. Halmos, \textsl{%
Lectures on Boolean Algebras}, Van Nostrand, 1963).\newline
\noindent \lbrack GofNN] C. Goffman, C. J. Neugebauer and T. Nishiura,
Density topology and approximate continuity, \textsl{Duke Math. J.} \textbf{%
28} (1961), 497--505.\newline
\noindent \lbrack GofW] C. Goffman, D. Waterman, Approximately continuous
transformations, \textsl{Proc. Amer. Math. Soc.} \textbf{12} (1961),
116--121.\newline
\noindent \lbrack GroE] K.-G. Grosse-Erdmann, An extension of the
Steinhaus-Weil theorem. \textsl{Colloq. Math.} \textbf{57}.2 (1989),
307--317.\newline
\noindent \lbrack Hal] P. R. Halmos, \textsl{Measure Theory}, Grad. Texts in
Math. \textbf{18}, Springer 1974. (1st. ed. Van Nostrand, 1950).\newline
\noindent \lbrack Has] H. Hashimoto, On the *topology and its application.
\textsl{Fund. Math. }\textbf{91} (1976), no. 1, 5--10.\newline
\noindent \lbrack HauP] O. Haupt, C. Pauc, La topologie approximative de
Denjoy envisag\'{e}e comme vraie topologie.\textsl{\ C. R. Acad. Sci. Paris}
\textbf{234} (1952), 390--392.\newline
\noindent \lbrack Hau] F. Hausdorff, \textsl{Grundz\"{u}ge der Mengenlehre},
Veit, Leipzig, 1914 (reprinted, Chelsea, N. York, 1949).\newline
\noindent \lbrack Hay] E. Hayashi, Topologies defined by local properties.
\textsl{Math. Ann.} \textbf{156} (1964), 205--215.\newline
\noindent \lbrack HewR] E. Hewitt, K. A. Ross, \textsl{Abstract Harmonic
Analysis,} Vol. I, Grundl. math. Wiss. \textbf{115}, Springer 1963 [Vol. II,
Grundl. 152, 1970].\newline
\noindent \lbrack HodHLS] W. Hodges, I. Hodkinson, D. Lascar, S. Shelah, The
small index property for $\omega $-stable $\omega $-categorical structures
and for the random graph, \textsl{J. London Math. Soc.} \textbf{48} (1993),
204--218.\newline
\noindent \lbrack HofT] J. Hoffmann-J\o rgensen, F. Tops\o e, Analytic
spaces and their application, in [Rog, Part 3].\newline
\noindent \lbrack Jab1] E. Jab\l o\'{n}ska, Some analogies between Haar
meager sets and Haar null sets in abelian Polish groups. \textsl{J. Math.
Anal. Appl.} \textbf{421} (2015),1479--1486.\newline
\noindent \lbrack Jab2] E. Jab\l o\'{n}ska, A theorem of Piccard's type in
abelian Polish groups.\textsl{\ Anal. Math.} \textbf{42} (2016), 159--164.%
\newline
\noindent \lbrack JanH] D. Jankovi\'{c}, T. R. Hamlett, New topologies from
old via ideals, \textsl{Amer. Math. Monthly} \textbf{97} (1990), 295-310.%
\newline
\noindent \lbrack KanK] R. Kannan, C.K. Krueger, \textsl{Advanced analysis
on the real line}, Universitext, Springer, 1996.\newline
\noindent \lbrack Kec] A. S. Kechris, \textsl{Classical descriptive set
theory}, Graduate Texts in Mathematics \textbf{156}, Springer, 1995.\newline
\noindent \lbrack Kem] J. H. B. Kemperman, A general functional equation,
\textsl{Trans. Amer. Math. Soc.} \textbf{86} (1957), 28--56.\newline
\noindent \lbrack Kha] A. B. Kharazishvili,\textsl{\ Transformation groups
and invariant measures. Set-theoretical aspects.} World Scientific, 1998.%
\newline
\noindent \lbrack Kom] Z. Kominek, On an equivalent form of a Steinhaus
theorem, \textsl{Math. (Cluj)} \textbf{30} (53)(1988), 25-27.\newline
\noindent \lbrack Kuc] M. Kuczma, \textsl{An introduction to the theory of
functional equations and inequalities. Cauchy's equation and Jensen's
inequality.} 2nd ed., Birkh\"{a}user, 2009 [1st ed. PWN, Warszawa, 1985].%
\newline
\noindent \lbrack KucS] M. Kuczma, J. Sm\'{\i}tal, On measures connected
with the Cauchy equation. \textsl{Aequationes Math.} \textbf{14} (1976), no.
3, 421--428.\newline
\noindent \lbrack LukMZ] J. Luke\v{s}, Jaroslav, J. Mal\'{y}, L. Zaj\'{\i}%
\v{c}ek, \textsl{Fine topology methods in real analysis and potential theory}%
. Lecture Notes in Math. \textbf{1189}, Springer, 1986.\newline
\noindent \lbrack Mar1] N. F. G. Martin, A topology for certain measure
spaces, \textsl{Trans. Amer. Math. Soc.} \textbf{112} (1964), 1-18.\newline
\noindent \lbrack Mar2] N. F. G. Martin, Generalized condensation points,
\textsl{Duke Math. J.}, \textbf{28} (1961), 507-514.\newline
\noindent \lbrack MarS] D. A. Martin, J. R. Steel, Projective determinacy.
\textsl{Proc. Nat. Acad. Sci. U.S.A.} \textbf{85} (1988), 6582--6586.\newline
\noindent \lbrack MatZ] E. Mato\u{u}skov\'{a}, M. Zelen\'{y}, A note on
intersections of non--Haar null sets, \textsl{Colloq. Math.} \textbf{96}
(2003), 1-4.\newline
\noindent \lbrack MilO] H. I. Miller and A. J. Ostaszewski, Group action and
shift-compactness, \textsl{J. Math. Anal. App.} \textbf{392} (2012), 23-39.%
\newline
\noindent \lbrack Mos] Y. V. Mospan, A converse to a theorem of Steinhaus,
\textsl{Real An. Exch.} \textbf{31} (2005), 291-294.\newline
\noindent \lbrack Mue] B. J. Mueller, Three results for locally compact
groups connected with the Haar measure density theorem, \textsl{Proc. Amer.
Math. Soc.} \textbf{16} (6) (1965), 1414-1416.\newline
\noindent \lbrack Mun] J. R. Munkres, \textsl{Topology, a first course},
Prentice-Hall, 1975.\newline
\noindent \lbrack Nik] O. Nikodym, Sur une propri\'{e}t\'{e} de l'op\'{e}%
ration $A,$ \textsl{Fund. Math.} \textbf{7} (1925), 149-154.\newline
\noindent \lbrack Ost1] A. J. Ostaszewski, Analytically heavy spaces:
Analytic Cantor and Analytic Baire Theorems, \textsl{Topology and its
Applications} \textbf{158} (2011), 253-275.\newline
\noindent \lbrack Ost2] A. J. Ostaszewski, Analytic Baire spaces, \textsl{%
Fund. Math.} \textbf{217} (2012), 189-210.\newline
\noindent \lbrack Ost3] A. J. Ostaszewski, Beyond Lebesgue and Baire III:
Steinhaus's Theorem and its descendants, \textsl{Topology and its
Applications} \textbf{160} (2013), 1144-1154.\newline
\noindent \lbrack Ost4] A. J. Ostaszewski, Almost completeness and the
Effros Theorem in normed groups, \textsl{Topology Proceedings} \textbf{41}
(2013), 99-110.\newline
\noindent \lbrack Ost5] A. J. Ostaszewski, Shift-compactness in almost
analytic submetrizable Baire groups and spaces, \textsl{Topology Proceedings}
\textbf{41} (2013), 123-151.\newline
\noindent \lbrack Ost6] A. J. Ostaszewski, Effros, Baire, Steinhaus and
non-separability, \textsl{Topology and its Applications (Mary Ellen Rudin
Memorial Volume)} \textbf{195} (2015), 265-274.\newline
\noindent \lbrack Oxt1] Oxtoby, John C. Invariant measures in groups which
are not locally compact. \textsl{Trans. Amer. Math. Soc.} \textbf{60},
(1946), 215--237.\newline
\noindent \lbrack Oxt2] J. C. Oxtoby: \textsl{Measure and category}, 2$^{%
\text{nd}}$ ed. Graduate Texts in Math. \textbf{2}, Springer, 1980 (1$^{%
\text{st}}$ ed. 1972).\newline
\noindent \lbrack Pat] A. L. T. Paterson, \textsl{Amenability}. Math.
Surveys and Mon. \textbf{29}, Amer. Math. Soc., 1988.\newline
\noindent \lbrack Rog] C. A. Rogers, J. Jayne, C. Dellacherie, F. Tops\o e,
J. Hoffmann-J\o rgensen, D. A. Martin, A. S. Kechris, A. H. Stone, \textsl{%
Analytic sets,} Academic Press, 1980.\newline
\noindent \lbrack Rud] W. Rudin, \textsl{Real and complex analysis}, 3$^{%
\text{rd}}$ ed. McGraw-Hill, 1987.\newline
\noindent \lbrack Sam] P. Samuels, A topology formed from a given topology
and ideal. \textsl{J. London Math. Soc.} (2) \textbf{10} (1975), 409--416.%
\newline
\noindent \lbrack Sch] S. Scheinberg, Topologies which generate a complete
measure algebra, \textsl{Adv. Math.} \textbf{7} (1971), 231-239. \newline
\noindent \lbrack ShiT] H. Shi, B. S. Thomson, Haar null sets in the space
of automorphisms on [0,1]. \textsl{Real Anal. Exchange} \textbf{24}
(1998/99), 337--350.\newline
\noindent \lbrack Sim] S. M. Simmons, A converse Steinhaus theorem for
locally compact groups, \textsl{Proc. Amer. Math. Soc.} \textbf{49 }(1975),
383-386.\newline
\noindent \lbrack Sol1] S. Solecki, Size of subsets of groups and Haar null
sets. \textsl{Geom. Funct. Anal.} \textbf{15} (2005), no. 1, 246--273.
\newline
\noindent \lbrack Sol2] S. Solecki, Amenability, free subgroups, and Haar
null sets in non-locally compact groups. \textsl{Proc. London Math. Soc.}
\textbf{(3) 93} (2006), 693--722.\newline
\noindent \lbrack Sol3] S. Solecki, A Fubini theorem. \textsl{Topology Appl.}
\textbf{154} (2007), no. 12, 2462--2464.\newline
\noindent \lbrack SteS] E. M. Stein and R. Shakarchi, \textsl{Real analysis,
Measure theory, integration and Hilbert spaces}, Princeton University press,
2005.\newline
\noindent \lbrack Str] K. Stromberg, An elementary proof of Steinhaus-Weil's
theorem, \textsl{Proc. Amer. Math. Soc.} \textbf{36} (6) (1972), 308.\newline
\noindent \lbrack TomW] G. Tomkowicz, S. Wagon, \textsl{The Banach-Tarski
paradox.} Cambridge University Press, 2016 (1$^{\text{st}}$ ed., S. Wagon,
1985).\newline
\noindent \lbrack Wei] A. Weil, \textsl{L'integration dans les groupes
topologiques}, Actualit\'{e}s Scientifiques et Industrielles 1145, Hermann,
1965.\newline
\noindent \lbrack Wil] W. Wilczy\'{n}ski, \textsl{Density topologies},
Handbook of measure theory, Vol. I, II, 675--702, North-Holland, 2002.%
\newline
\noindent \lbrack WilK] W. Wilczy\'{n}ski, A. B. Kharazishvili, Translations
of measurable sets and of sets having the Baire property. (Russian) \textsl{%
Soobshch. Akad. Nauk Gruzii} \textbf{145.1} (1992), no. 1, 43--46.\newline
\noindent \lbrack Woo] W. H. Woodin,\textsl{\ The axiom of determinacy,
forcing axioms, and the nonstationary ideal.} 2$^{\text{nd}}$ ed., de
Gruyter Series in Logic and its Applications\textbf{\ 1}. de Gruyter, 2010.%
\newline

\bigskip

\noindent Mathematics Department, Imperial College, London SW7 2AZ;
n.bingham@ic.ac.uk \newline
Mathematics Department, London School of Economics, Houghton Street, London
WC2A 2AE; A.J.Ostaszewski@lse.ac.uk

\bigskip

\newpage

\section*{\qquad \qquad \qquad \qquad Appendix\protect\footnote{%
This Appendix for the arXiv version only.}\newline
\textbf{The Steinhaus Property }$AA^{-1}$\textbf{\ versus the Steinhaus
Property }$AB^{-1}$}

In this section we expand on the comment in \S 6.9 and further clarify the
relation between two versions of the Steinhaus interior points property: the
simple version concerning sets $AA^{-1}$ and the composite, more embracing
one, concerning sets $AB^{-1}$, for sets from a given family $\mathcal{H}$.
The latter is connected to a strong form of metric transitivity: Kominek
[Kom] shows for a general separable Baire topological group $G$ equipped
with an inner-regular measure $\mu $ defined on some $\sigma $-algebra $%
\mathcal{M}$ that $AB^{-1}$ has non-empty interior for all $A,B\in \mathcal{M%
}_{+}$ iff for each countable dense set $D$ and each $E\in \mathcal{M}_{+}$
the set $X\backslash DE\in \mathcal{M}_{0}$. Care is required when moving to
the alternative property $AB,$ since the family $\mathcal{H}$ need not be
preserved under inversion. In general, the simple property does not imply
the composite,\ as noted in in \S 6.9: Mato\u{u}skov\'{a} and Zelen\'{y}
[MatZ] show that in any non-locally compact abelian Polish group there are
closed non-(left) Haar null sets $A,B$ such that $A+B$ has empty interior;
Jab\l o\'{n}ska [Jab2] has shown that likewise in any non-locally compact
abelian Polish group there are closed non-Haar meager sets $A,B$ such that $%
A+B$ has empty interior.

Below we identify some conditions on a family of sets $A$ with the simple $%
AA^{-1}$ property which do indeed imply the $AB^{-1}$ property. What follows
is a generalization to a group context of relevant observations from [BinO5]
from the classical context of $\mathbb{R}$.

The motivation for the definition below is that the space $H\ $is a subgroup
of a topological group $G$ from which it inherits a (necessarily)
translation-invariant (i.e. on either sided) topology $\tau .$ Various
notions of `density at a point' give rise to `density topologies' as above,
which are translation-invariant since they may be obtained via translation
to a fixed reference point: early examples, which originate in spirit with
Denjoy as interpreted by Haupt and Pauc [HauP], were studied intensively in
[GofW], [GofNN], soon followed by [Mar1] and [Mue]; more recent examples
include [FilW] and others investigated by the Wilczy\'{n}ski school, cf.
[Wil].

Proposition A below embraces as an immediate corollary the case $H=G$ with $%
G\ $locally compact and $\sigma $ the Haar density topology (see [BinO4]).
Proposition B proves that Proposition A applies also to the ideal topology
(in the sense of [LukMZ]) generated from the ideal of Haar null sets of an
abelian Polish group.

We recall that a group $H$ carries a \textit{left semi-topological}
structure $\tau $ if the topology $\tau $ is left-translation invariant
[ArhT] ($hU\in \tau $ iff $U\in \tau $); the structure is \textit{%
semi-topological} if it is also right-invariant, i.e. briefly: $\tau $ is
translation invariant. $H$ is a \textit{quasi-topological group} under $\tau
$ if $\tau $ is both left and right invariant and inversion is $\tau $%
-continuous.

\bigskip

\noindent \textbf{Definition.} For $H$ a group with a translation-invariant
topology $\tau ,$ call a topology $\sigma \supseteq \tau $ a \textit{%
Steinhaus refinement }if:

\noindent i) \textrm{int}$_{\tau }$($AA^{-1})\neq \emptyset $ for each
non-empty $A\in \sigma ,$ and

\noindent ii) $\sigma $ is involutive-translation invariant: $hA^{-1}\in
\sigma $ for all $A\in \sigma $ and all $h\in H.$

\bigskip

Property (ii) above (called simply `invariance' in [BarFN]) apparently calls
for only left invariance, but in fact, via double inversion, delivers
translation invariance, since $Uh=(h^{-1}U^{-1})^{-1}$; then $H$ under $%
\sigma $ is a semi-topological group with a continuous inverse, so a \textit{%
quasi-topological group}.

\bigskip

\noindent \textbf{Proposition A.} \textit{If }$\tau $\textit{\ is
translation-invariant, and }$\sigma \supseteq \tau $\textit{\ is a Steinhaus
refinement topology, then }\textrm{int}$_{\tau }$($AB^{-1})\neq \emptyset $%
\textit{\ for non-empty }$A,B\in \sigma .$ \textit{In particular, as }$%
\sigma $\textit{\ is preserved under inversion, also }\textrm{int}$_{\tau }$(%
$AB)\neq \emptyset $\textit{.for }$A,B\in \sigma $.\textit{\ }

\bigskip

\noindent \textbf{Proof.} Suppose $A,B\in \sigma $ are non-empty; as $%
B^{-1}\in \sigma ,$ choose $a\in A$ and $b\in B,$ then by (ii)%
\[
1_{H}\in C:=a^{-1}A\cap b^{-1}B^{-1}\in \sigma .
\]%
By (i), for some non-empty $W\in \tau ,$%
\[
W\subseteq CC^{-1}=(a^{-1}A\cap b^{-1}B)\cdot (A^{-1}a\cap B^{-1}b)\subseteq
(a^{-1}A)\cdot (B^{-1}b).
\]%
As $\tau $ is translation invariant, $aWb^{-1}\in \tau $ and
\[
aWb^{-1}\subseteq AB^{-1},
\]%
the latter since for each $w\in W$ there are $x\in A,y\in B^{-1}$ with
\[
w=a^{-1}x.yb:\qquad awb^{-1}=xy\in AB^{-1}.
\]%
So \textrm{int}$_{\tau }$($AB^{-1})\neq \emptyset .$ $\square $

\bigskip

\noindent \textbf{Corollary.} \textit{In a locally compact group the Haar
density topology is a Steinhaus refinement.}

\bigskip

\noindent \textbf{Proof. }Property (i) follows from Weil's theorem since
density open sets are non-null measurable; left translation invariance in
(ii) follows from left invariance of Haar measure, while involutive
invariance holds, as any measurable set of positive Haar measure has
non-null inverse. $\square $

\bigskip

A weaker version, inspired by metric transitivity, comes from applying the
following concept.

\bigskip

\noindent \textbf{Definition.} Say $H$ acts transitively on $\mathcal{H}$
for each $A,B\in \mathcal{H}$ if there is $h\in H$ with $A\cap hB\in
\mathcal{H}$.

\bigskip

Thus a locally compact topological group acts transitively on the non-null
Haar measurable sets (in fact, either sidedly); this follows from Fubini's
Theorem [Hal, 36C], via the average theorem [Hal, 59.F]:%
\[
\int_{G}|g^{-1}A\cap B|dg=|A|\cdot |B^{-1}|\qquad (A,B\in \mathcal{M}),
\]%
($g=ab^{-1}$ iff $g^{-1}a=b)$ -- cf. [TomW, \S 11.3 after Th. 11.17].

[MatZ] show that in any non-locally compact abelian Polish group $G$ there
exist two non-Haar null sets, $A,B\notin \mathcal{HN}$, such that $A\cap
hB\in \mathcal{HN}$ for all $h;$ that is, $G$ does not act transitively on
the non-Haar null sets.

\bigskip

\noindent \textbf{Proposition A}$^{\prime }$\textbf{.} \textit{In a group }$%
(H,\tau )$ \textit{with }$\tau $\textit{\ translation-invariant, if }$H$
\textit{acts transitively on a family of subsets }$\mathcal{H}$\textit{\
with the simple Steinhaus property, then }$\mathcal{H}$ \textit{has the
composite Steinhaus property: }\textrm{int}$_{\tau }$($AB^{-1})\neq
\emptyset $\textit{\ for }$A,B\in \mathcal{H}$.\textit{\ Furthermore, if }$%
\mathcal{H}$\textit{\ is preserved under inversion, then also }\textrm{int}$%
_{\tau }$($AB)\neq \emptyset $\textit{\ for }$A,B\in \mathcal{H}$.

\bigskip

\noindent \textbf{Proof. }Argue as in Prop. A, but now for $A,B\in \mathcal{H%
}$ choose $h$ with $C:=A\cap hB\in \mathcal{H}$; then%
\[
CC^{-1}h=(A\cap hB)(A^{-1}\cap B^{-1}h^{-1})\subseteq AB^{-1}.\qquad \square
\]

\bigskip

\noindent \textbf{Definition.} In a quasi-topological group $(H,\tau )$ say
that a proper $\sigma $-ideal $\mathcal{H}$ has the \textit{Simple Steinhaus
Property }$AA^{-1}$ if $AA^{-1}$ has interior points for universally
measurable subsets $A\notin \mathcal{H}.$ This follows [BarFN].

\bigskip

\noindent \textbf{Proposition B.} \textit{If }$(H,\tau )$\textit{\ is a
quasi-topological group (i.e. }$\tau $\textit{\ is invariant with continuous
inversion) carrying a left invariant }$\sigma $\textit{-ideal }$\mathcal{H}$
\textit{with the Steinhaus property and }$\tau \cap \mathcal{H}=\{\emptyset
\}$\textit{, then the ideal-topology }$\sigma $\textit{\ with basis }%
\[
\mathcal{B}:=\{U\backslash N:U\in \tau ,N\in \mathcal{H}\}
\]%
\textit{is a Steinhaus refinement of }$\tau $\textit{.}

\textit{In particular, for }$(H,\tau )$\textit{\ an abelian Polish group,
the ideal topology generated by the }$\sigma $\textit{-ideal of Haar null
subsets is a Steinhaus refinement.}

\bigskip

\noindent \textbf{Proof. }If $U,V\in \mathcal{B}$ and $w\in U\cap V,$ choose
$M,N\in \mathcal{H}$ and $W_{M},W_{N}\in \tau $ such that $x\in
(W_{M}\backslash M)\subseteq U$ and $x\in (W_{N}\backslash N)\subseteq V,$
then as $M\cup N\in \mathcal{H}$
\[
x\in (W_{M}\cap W_{N})\backslash (M\cup N)\in \mathcal{B}.
\]%
So $\mathcal{B}$ generates a topology $\sigma $ refining $\tau .$ With the
same notation, $hU=hW_{M}\backslash hM\in \sigma ,$ as $hM\in H,$ and $%
U^{-1}=W_{M}^{-1}\backslash M^{-1}$; finally $UU^{-1}$ has non-empty $\tau $%
-interior, as $U\notin \mathcal{H}$ and is non-empty.

As for the final assertion concerned with an abelian Polish group context,
note that if $N$ is Haar null ($N\in \mathcal{HN}$), then $\mu (hN)=0$ for
some probability measure $\mu $ and all $h\in H,$ so $hN\in \mathcal{HN}$
for all $h\in H;$ furthermore, if $A\notin \mathcal{HN}$ then $A^{-1}\notin
\mathcal{HN}$, otherwise $\mu (hA^{-1})=0$ for some probability measure $\mu
$ and all $h\in H,$ and then, taking $\bar{\mu}(B)=\mu (B^{-1})$ for Borel $%
B,$ we have $\bar{\mu}(A)=0$ and $\bar{\mu}(hA)=\mu (A^{-1}h^{-1})=0$ for
all $h\in H,$ a contradiction. $\square $

\bigskip

\noindent \textbf{Remark.} A left Haar null set need not be right Haar null:
for one example see [ShiT], and for more general non-coincidence see Solecki
[Sol1, Cor. 6]. So the argument in Prop. B does not extend to the family of
left Haar null sets $\mathcal{LHN}$ of a \textit{non-commutative} Polish
group. Indeed, Solecki [Sol2, Th. 1.4] shows in the context of a countable
product of countable groups that the simpler Steinhaus property holds for $%
\mathcal{HN}$ iff $\mathcal{HN}=\mathcal{LHN}$.

\bigskip

We close with a result from [Kom]. Recall that $\mu $ is \textit{%
quasi-invariant} if $\mu $-nullity is translation invariant. The
transitivity assumption (of co-nullity) is motivated by \textit{Sm\'{\i}%
tal's lemma}, hich refers to a countable dense set -- see [KucS].

\bigskip

\noindent \textbf{Theorem K (}[Kom, Th. 5]).\textit{\ If }$\mu \in \mathcal{P%
}(G)$\textit{\ is quasi-invariant and there exists a countable subset }$%
H\subseteq G\ $\textit{with }$HM$\textit{\ co-null for all }$M\in \mathcal{M}%
_{+}(\mu ),$ \textit{then} \textrm{int}$(AB^{-1})\neq \emptyset $ \textit{%
for all }$A,B\in \mathcal{M}_{+}(\mu ).$

\bigskip

\noindent \textbf{Proof. }By regularity we may assume $A,B\in \mathcal{M}%
_{+}(\mu )$ are compact, so $AB^{-1}$ is compact. Fix $g\in G$; then by
quasi-invariance $\mu (gB)>0,$ so by the transitivity assumption, both $%
G\backslash HgB$ and $G\backslash HA$ are null, and so $HA\cap HgB\neq
\emptyset .$ Say $h_{1}a=h_{2}gb,$ for some $a\in A,b\in B,h_{1},h_{2}\in H,$
then $g=h_{2}^{-1}h_{1}ab^{-1}.$ As $g$ was arbitrary,%
\[
G=\dbigcup\nolimits_{h\in H}h_{2}^{-1}h_{1}AB^{-1}.
\]%
By Baire's Theorem, as $H$ is countable, \textrm{int}$(AB^{-1})\neq
\emptyset .$ $\square $

\end{document}